\documentstyle[12pt]{article}
\input epsf
\newcounter{THNO}[section]
\renewcommand{\theTHNO}{\arabic{section}.\arabic{THNO} }

\def\th#1{\refstepcounter{THNO}\par\vspace{0.7cm}
\par\noindent\begingroup \it
\leftskip=0em\hspace{0em}{\bf\ #1 \theTHNO\ }}
\def\eth{\par\endgroup}
\def\pr{\par\noindent{\bf\ Proof. }}
\def\db#1{ D^b_{coh}({#1})}
\def\h#1,#2{{\rm Hom}({#1}\:,\; {#2})}
\def\H#1,#2,#3,#4{{\rm Hom}^{#1}_{#2}({#3}\:,\; {#4})}
\def\lto{\longrightarrow}
\def\o#1{{\cal O}_{#1}}
\def\b{\; \epsffile{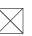}\;}

\begin{document}
\begin{center}
{\LARGE Semiorthogonal decompositions}\\

\vspace{4ex}

{\LARGE for algebraic varieties.}
\end{center}

\vspace{2ex}

\centerline{
\begin{tabular}{ll}
{\large A.~Bondal}\qquad&\qquad{\large D.~Orlov}\\[7pt]
Algebra section,&\qquad Slaviansky br.\\
Steklov Math Institute\qquad&\qquad d.~9, corp. 4\\
Vavilova 42,&\qquad kv.80 \\
117333 Moscow&\qquad 121352 Moscow\\
RUSSIA&\qquad RUSSIA
\end{tabular}}

\vspace{5cm}

\begin{center}
{\large  March 6, 1995}
\end{center}
\newpage
\section*{ SEMIORTHOGONAL DECOMPOSITIONS FOR ALGEBRAIC VARIETIES}
\begin{abstract}
A criterion for a functor between derived categories of coherent sheaves to be full and faithful is given. A semiorthogonal decomposition for the derived category of coherent sheaves on the intersection of two even dimensional quadrics is obtained. The behaviour of derived categories with respect to birational transformations is investigated.
A theorem about reconstruction of a variety from the derived category of coherent sheaves is proved.
\end{abstract}

\tableofcontents


\section{Introduction.}

This paper is devoted to study of the derived categories of coherent sheaves on smooth algebraic varieties. Of special interest for us is the case when there exists a functor ${\db M}\lto {\db X}$ which is full and faithful.

It appears that some geometrically important constructions for moduli spaces of (semistable) coherent sheaves on varieties can be interpreted as instances of this situation. Conversely, we are convinced that any example of such a functor is geometrically meaningful.

If a functor $\Phi: {\db M}\lto {\db X}$ is full and faithful, then it induces  a semiorthogonal decomposition (see definition in ch.2) of ${\db X}$ with the 2--step chain $\Bigl( {\db M}^{\perp}$, ${\db M} \Bigl)$, where ${\db M}^{\perp}$ is the right orthogonal to ${\db M}$ in ${\db X}$.

Decomposing summands of this chain, one can obtain a  semiorthogonal decomposition with arbitrary number of steps.
Full exceptional sequences existing on some Fano varieties (see \cite{K}) provide with examples of such decompositions. For this case, every step of the chain is equivalent to the derived category of vector spaces or, in other words, sheaves over the point.

This leads to the idea that the derived category of coherent sheaves might be reasonable to consider as an incarnation of the motive of a variety, while semiorthogonal decompositions are a tool for  simplification of a motive similar to spliting by projectors in the Grothendieck motive theory.

Main result of ch.1 is a criterion for fully faithfulness. Roughly speaking, it claims that for a functor ${\db M}\lto{\db X}$ to be full and faithful it is sufficient to satisfy this property on the full subcategory of the one dimensional skyscraper sheaves and its translations.

Let us mention that ${\db M}^{\perp}$ might be zero. In this case we obtain an equivalence of derived categories ${\db M}\stackrel{\sim}{\lto}{\db X}$. Examples of such equivalences have been considered by Mukai in \cite{Mu}, \cite{Mu2} (see ch.1). In ch.3 we prove such equivalences for some flop birational transformations.

Ch.2 is devoted to description of a semiorthogonal decomposition for ${\db X}$, when $X$ is the smooth intersection of two even dimensional quadrics. It appears that if we consider the hyperelliptic curve $C$ which is a double covering of the projective line parametrizing the pencil of quadrics, with ramification in the points corresponding to degenerate quadrics, then ${\db C}$ is embedded in ${\db X}$ as a full subcategory. The orthogonal to ${\db C}$ in ${\db X}$ is decomposed in an exceptional sequence (of linear bundles ). This allows to identify moduli spaces of semistable bundles (of arbitrary rank) on the curve with moduli spaces of complexes of coherent sheaves on the intersection of quadrics. For rank 2 bundles such identification is well known (see \cite{DR}) and was used for computation of cohomologies of moduli spaces \cite{Bar} and for verification of the Verlinde formula.

In ch.3 we investigate the behaviour of ${\db X}$ under birational transformations. We prove that for a couple of varieties $X$ and $X^+$ related by some flips the category ${\db {X^+}}$ has a natural full and faithful embedding in ${\db X}$. This suggests the idea that the minimal model program of the birational geometry can be considered as a `minimization' for the derived category of coherent sheaves in a given birational class. 

 We also explore some cases of flops. Considered examples allow us to state a conjecture that the derived categories of coherent sheaves on varieties connected by a flop are equivalent.

Examples of varieties having equivalent derived categories appeal to the question:

{\it to which extent a variety is determined by its derived category?}

In ch.4 we prove a reconstruction theorem, which claims that
if $X$ is a smooth algebraic variety with ample either canonical or anticanonical sheaf, then another algebraic variety $X'$ having equivalent the derived category of coherent sheaves ${\db X}\simeq {\db {X'}}$ should be biregulary isomorphic to $X$.

As a by-product we obtain a description for the group of
 auto-equivalences of ${\db X}$ provided $X$ has ample either canonical or anticanonical class.

We are grateful to Max--Planck--Institute for hospitality and stimulating atmosphere. Our special thanks go to 
S.Kuleshov for  help during preparation of this paper. The work was partially supported by International Science Foundation Grant M3E000 and Russian Fundamental Research Grant.

{\section {Full and faithful functors.}}

\hspace*{0.6cm}For a smooth algebraic variety  $X$ over an  algebraically closed field $k$ of characteristic $0$ by
 ${\db X}$ (resp., $D^b_{Qcoh}(X)$) we denote the bounded derived category of
coherent (resp., quasicoherent) sheaves over $X$. Notations like $f^{*}, f_{*},
\otimes, {\rm Hom}, {\cal H}om$ etc. are reserved for derived functors between derived categories, whereas $R^if_*, {\rm Hom}^i$, etc. (resp., $L^if^*$) denote i--th (resp., (-i)--th) cohomology of a complex obtained by applying $f_*, {\rm Hom}$ etc. (resp., $f^*$); $[n]$ denotes the translation by $n$ functor in a triangulated category.

Let $X$ and $M$ be smooth algebraic varieties of dimension $n$ and $m$ respectively, and $E$ an object of ${\db {X \times M}}$. With $E$ one can associate a couple of functors 
$$
 \Phi_{E} : {\db M}\lto {\db X},
$$
$$
\Psi_{E} : {\db X}\lto {\db M}.
$$
Denote by $p$ and $\pi$ the projections of $M \times X$ to $M$ and $X$ respectively.
$$
\begin{array}{ccc}
M\times X&\stackrel{\pi }{\longrightarrow}& X\\
\llap{\footnotesize $p$} \downarrow &&\\
M &&
\end{array}
$$

Then $\Phi_{E}$ and $\Psi_{E}$ are defined by the formulas:
$$
\Phi_{E}(\cdot):=\pi_* (E \otimes p^*(\cdot)),
$$
$$
\Psi_{E}(\cdot):=p_* (E \otimes \pi^*(\cdot)).
$$
The main goal of this chapter is the proof of the following

\th{Theorem}\label{mai} 
Let $M$ and $X$ be smooth algebraic varieties and\hfill\\ 
 $E\in{\db {M\times X}}$. Then $\Phi_{E}$ is  full and faithful functor,
if and only if the following orthogonality conditions are verified: 
$$
\begin{array}{lll}
i) & {\H i, X, \Phi_E({\cal O}_{t_1}), {\Phi_{E}({\cal O}_{t_2})}} = 0 & \qquad \mbox{for every }\: i\;\mbox{ and } t_1\ne t_2.\\ 
&&\\
ii) & {\H 0, X, \Phi_E({\cal O}_t), {\Phi_E({\cal O}_t)}} = k,&\\ 
&&\\
& {\H i, X, \Phi_E({\cal O}_t), {\Phi_E({\cal O}_t)}} = 0 , & \qquad \mbox{ for }i\notin [0, dim M].
\end{array}
$$
Here $t$, $t_{1}$, $t_{2}$ are points of $M$, ${\cal O}_{t_{i}}$  corresponding skyscraper sheaves.
\eth
\bigskip

Let us mention that if some full subcategory ${\cal C}\subset {\cal D}$ generates ${\cal D}$ as a triangulated category then for an exact functor ${\cal D}\lto{\cal D}'$ to be full and faithful it is sufficient to be full on ${\cal C}$. Unfortunately, the class of skyscraper sheaves does not generate $\db{M}$ as a triangulated category if $dimM>0$. At the level of the Grothendieck group $K_0(M)$ they generate only the lowest term of the topological filtration.
  
The proof of the theorem is preceded by a series of assertions concerning functors between and objects from the derived categories of complexes of coherent sheaves on smooth varieties.

For any object $E$ from $\db{X}$ we denote by $E^{\vee}$ the dual object:
$$
E^{\vee}:={\cal H}om( E,\:{\o X}).
$$


\th{Lemma}\label{adj}
The left adjoint functor to $\Phi_{E}$ is
$$
\Psi_{E^{\vee}\otimes \pi^*\omega_X}[n]:= 
  p_{*}(E^{\vee}\otimes \pi^*\omega_X\otimes \pi^*(\cdot))[n]. 
$$
\eth
{\bf Proof} is given by a series of natural isomorphisms, which come from the adjoint property of functors and Serre duality:
$$
\begin{array}{l}
 {\h A, {\pi_*(E\otimes p^*B)}}\cong \\ 
{\h {\pi^*A}, {E\otimes p^* B}}\cong  \\
{\h {p^* B}, {\pi^* A\otimes E^{\vee}\otimes \omega_{X\times M}[n+m]}}^* \cong \\
{\h B, {p_*(\pi^*(A\otimes \omega_X[n])\otimes E^{\vee})\otimes \omega_M[m]}}^* \cong\\
 {\h {p_*(\pi^*(A\otimes \omega_X[n])\otimes E^{\vee})}, B}.
\end{array}
$$
\bigskip

The next lemma differs from analogous in \cite{H} in what concerns base change (we consider arbitrary $g$ instead of flat one in \cite{H}) and morphism $f$ (we consider only smooth morphism instead of arbitrary one in \cite{H}).


\th{Lemma}\label{isof}
Let $f: X\to Y $ be a smooth morphism of relative dimension $r$ of smooth projective varieties and $ g:Y'\to Y $  a base change, with
$Y^{\prime}$ being a smooth variety. Define $X'$ as the cartesian product $X'=X\times_Y Y'$.
$$
\begin{array}{ccc}
X^{\prime}=X\times_Y Y^{\prime}& \stackrel{ g'}{\longrightarrow}& X\\      
\llap{$f'$}\downarrow &&
\llap{$f$}\downarrow \\
Y'&\stackrel{g}{\longrightarrow}& Y\\
\end{array}
$$

Then there is a natural isomorphism of functors:
$$
g^*f_*(\cdot) \simeq {f^{\prime}}_* {g^{\prime}}^*(\cdot).
$$
\eth
\pr
First, note that the right adjoint functors to $g^*f_*$ and $ f'_* g'^*$ are, respectively, $f^{!}g_*$ and $g'_*f'^!$ , where $f^!$ denote the right adjoint functor to $f_*$. We are going to prove that $f^{!}g_*$ and $g'_*f'^!$ are isomorphic.

Serre duality gives a natural isomorphism
\begin{equation}\label{1}
f^!(\cdot)\simeq f^*(\cdot)\otimes \omega_{X/Y}[r].
\end{equation}
Hence,
\begin{equation}\label{2}
f^!g_*(\cdot)\simeq f^*g_*(\cdot)\otimes \omega_{X/Y}[r].
\end{equation}
Analogously,
$$
g'_*f'^!(\cdot)\simeq g'_*(f'^*(\cdot)\otimes \omega_{X'/Y'}[r])\simeq g'_*(f'^*(\cdot)\otimes g'^*\omega_{X/Y}[r]).
$$
The latter isomorphism goes from the fact that for a smooth $f$ differentials are compatible with base change (see \cite{H},III,\S1,p.141).
Then, by the projection formula one has
\begin{equation}\label{3}
g'_* f'^!(\cdot)\simeq g'_* f'^*(\cdot)\otimes \omega_{X/Y}[r].
\end{equation}
By the theorem of flat base change (see \cite{H},II,\S5,prop.5.12) one has
$$
g'_* f'^*\simeq f^*g_*.
$$
Formulas (\ref{2}) and (\ref{3}) imply a functorial isomorphism of $g'_* f'^!(\cdot)$ and $f^! g_*(\cdot)$.
Therefore, $g^*f_*(\cdot)$ is isomorphic to $f'_* g'^*(\cdot)
$.
\bigskip

Let $X, Y, Z$ be smooth projective varieties and $I, J, K$ objects of $\db{X\times Y}$, $\db{Y\times Z}$ and $\db{X\times Z},$ respectively. Consider the following diagram of projections

\hspace*{5cm}\epsffile{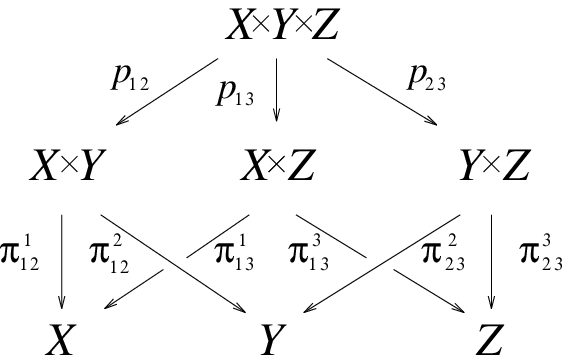}

and the triple of functors
$$
 \phi_{I} : {\db X}\lto {\db Y},
$$
$$
 \psi_{J} : {\db Y}\lto {\db Z},
$$
$$
 \chi_{K} : {\db X}\lto {\db Z},
$$
defined by the formulas
$$
\phi_I={\pi_{12}^2}_{*}(I\otimes {\pi_{12}^{1}}^{*}(\cdot)),
$$
$$
\psi_J={\pi^3_{23}}_*(J\otimes {\pi^2_{23}}^*(\cdot)),
$$
$$
\chi_K={\pi^3_{13}}_*(K\otimes {\pi^1_{13}}^*(\cdot)).
$$

The next proposition from \cite{Mu} is an analog for derived categories of the composition law for correspondences (see \cite{M}).

\th {Proposition}\label{comp}
The composition functor for $\phi_I$ and $\psi_J$ is isomorphic to $\chi_K$ with
$$
K={p_{13}}_*({p_{23}}^*J\otimes {p_{12}}^*I).
$$ 
\eth
\pr
It goes from the following sequence of natural isomorphisms, which uses the projection formula and a base change theorem from \cite{H}:
$$
\begin{array}{l}
\psi_J \circ \phi_I(\cdot)\cong {\pi^3_{23}}_*(J\otimes {\pi^2_{23}}^*({\pi^2_{12}}_*(I\otimes {\pi^1_{12}}^*(\cdot))))
\cong \\
 {\pi^3_{23}}_*(J\otimes {p_{23}}_*({p_{12}}^*(I\otimes {\pi^1_{12}}^*(\cdot))))\cong\\
 {\pi^3_{23}}_*{p_{23}}_*({p_{23}}^*J\otimes {p_{12}}^*(I\otimes {\pi^1_{12}}^*(\cdot)))\cong\\
 {\pi^3_{13}}_*{p_{13}}_*({p_{23}}^*J\otimes {p_{12}}^*I\otimes {p_{12}}^* {\pi^1_{12}}^*(\cdot))\cong\\
 {\pi^3_{13}}_*{p_{13}}_*({p_{23}}^*J\otimes {p_{12}}^*I\otimes {p_{13}}^* {\pi^1_{13}}^*(\cdot))\cong\\
 {\pi^3_{13}}_*({p_{13}}_*({p_{23}}^*J\otimes {p_{12}}^*I)\otimes {\pi^1_{13}}^*(\cdot)).\\
\end{array}
$$

\th{Proposition}\label{tor}
Let $j: Y\hookrightarrow X$ be a smooth irreducible subvariety of codimension $d$ of a smooth algebraic variety $X$, and $K$ a non-zero object of $\db{X}$ satisfying following conditions:

a) $i^*_x K=0$, \qquad for any closed point $x\stackrel{i_x}{\hookrightarrow} X\setminus Y$, 

b) $L^i i^*_x K=0$,\qquad when $ i\notin [0, d]$, for any closed point $x\stackrel{i_x}{\hookrightarrow} Y$. 

Then

i) $K$ is a pure sheaf (i.e. quasiisomorphic to its zero cohomology sheaf),

ii) the support of $K$ is $Y$.
\eth
\pr
Let ${\cal H}^q$ be the q--th cohomology sheaf of $K$. Then, for any point $x\stackrel{i_x}{\hookrightarrow} X$ there is spectral sequence with the $E_2$--term consisting of $L^p i^*_x({\cal H}^q)$ and converging to cohomology sheaves of $i_*(K)$:
$$
E^{-p,q}_2=L^p i^*_x({\cal H}^q) \Rightarrow L^{p-q} i^*_x(K)
$$
Recall that $L^i f^*$ denotes the (--i)--th cohomology of $f^*$ in accordance with notations of the analogous left derived functors between abelian categories.

If ${\cal H}^{q_{max}}$ is a non--zero sheaf with maximal $q$, then $L^0 i^*_x{\cal H}^{q_{max}}$ is intact by differentials while going to $E_{\infty}$. By assumptions of the proposition $L^q i^*_x K=0$, for $q>0$ and for any point $x\in X$. This implies $q_{max}\le 0$.

Considering the sheaf ${\cal H}^q$ with maximal $q$ among those having the support outside $Y$, one obtains by the same reasoning that all ${\cal H}^q$ actually have their support in $Y$.

Let ${\cal H}^{q_{min}}$ be the non--zero sheaf with minimal $q$. The spectral sequence is depicted in the following diagram:

\begin{figure}[th]
\hspace*{5cm}\epsffile{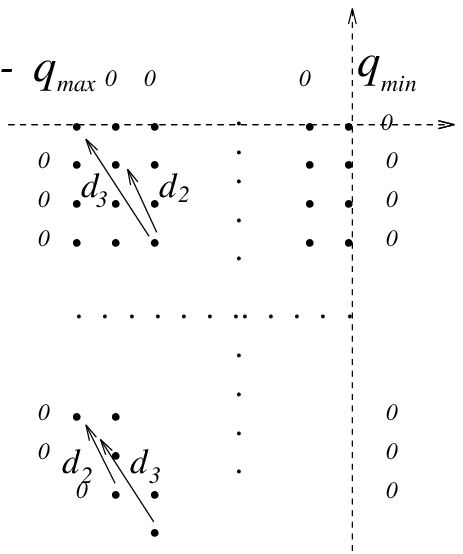}
\end{figure}

Consider any component $C\subset Y$ of the support of ${\cal H}^{q_{min}}$. If $c$ is the codimension of $ C$ in $X$, then $L^c i^*_{x_0}({\cal H}^{q_{min}})\ne 0$ for a general closed point $x_0\in C$.
It could have been killed in the spectral sequence only by $L^p i^*_x({\cal H}^q)$ with $p\ge c+2$. But for any sheaf $F$ the closed subscheme $S_m(F)$ of points of cohomological dimension $\ge m$ (see \cite{S})
$$
S_m(F)=\Bigl\{ x\in X\;\Bigl|\;L^p i^*_x(F)\ne 0,\quad\mbox{ for some}\:p\ge m\Bigl\}
$$
has codimension $\ge m$. Therefore, $S_m({\cal H})$ with $m\ge c+2$ cannot cover $C$, i.e. there exists a point $x_0\in C$, such that $L^c i^*_{x_0}({\cal H}^{q_{min}})$ survives at infinity in the spectral sequence, hence $L^{c-q_{min}} i^*_{x_0}(K)\ne 0$.

Then, by assumption b) of the proposition it follows that $c-q_{min}\le d$. Since $C$ belongs to $Y$, $c\ge d$, hence $q_{min}\ge 0$. In other words, $q_{min}=q_{max}$ and $K$ has the only non--trivial cohomology sheaf ${\cal H}^0$. This proves i).

Now consider ${\cal L}^i=L^i j^*K$. There is a spectral sequence for composition of $i^*_x$ and $j^*$:
$$
E^{-p,-q}_2=L^p i^*_x({\cal L}^q) \Rightarrow L^{p+q} i^*_x(K).
$$

Let ${\cal L}^{q_0}$ be a non--zero sheaf with maximal $q$. Since the support of $K$ belongs to $Y$, $q_0\ge d$. Again consider a component of the support for ${\cal L}^{q_0}$.
The same reasoning as above shows that if this component is of codimension $b$, then for some point $x_0$ in it, $L^b i^*_{x_0}({\cal L}^{q_0})$ survives in $E_{\infty}$ of the latter spectral sequence. By the assumptions of the proposition we have $q_0+b\le d$. This implies
$q_0=d$ and $b=0$. This means that the support of ${\cal L}^d$ is the whole $Y$. It follows that the support of $K$ coincides with $Y$. The proposition is proved.
\bigskip


{\bf Proof of the Theorem \ref{mai}.}
First, let us mention that if $\Phi_E$ is full and faithful functor, then conditions i) and ii) are verified for obvious reasons. Indeed, it is well known fact that extension groups between skyscraper sheaves in $\db{M}$ have the following form:
$$
\begin{array}{lll}
i) & {\H i, X, {\cal O}_{t_1}, {{\cal O}_{t_2}}} = 0 & \qquad \mbox{for every }\: i\;\mbox{ and } t_1\ne t_2;\\ 
ii) & {\H i, X, {\cal O}_t , {{\cal O}_t}} = \Lambda^i T_{M,t} , & \qquad \mbox{ for }i\in [0, dim M],\\
 & {\H i, X, {\cal O}_t , {{\cal O}_t}} = 0 , & \qquad \mbox{ for }i\notin [0, dim M].
\end{array}
$$
Here $t, t_1, t_2$ are points of $M$, $T_{M,t}$ the tangent vector space to $M$ at $t$, and $ \Lambda^i$ the $i$--th exterior power.

Fully faithfulness of $\Phi_E$ implies that the same relations are valid for images $\Phi_E({\cal O}_t)$ in $\db{X}$.

In what follows we prove the inverse statement.
 
 Consider composition of $\Phi_E$ with its left adjoint functor $\Phi_E^*$. 
We are going to prove that the canonical natural transformation
$\alpha : \Phi_E^*\circ\Phi_E\to id$ is an isomorphism of functors. This is equivalent to fully faithfulness of $\Phi_E$. Indeed, for any pair of objects $A, B\in\db{M}$ the natural homomorphism
$$
{\h A, B}\lto {\h \Phi_E A, {\Phi_E B}}\cong{\h \Phi_E^*\Phi_E A, B},
$$
is induced by $\alpha$.

 By lemma \ref{adj} we have
$$
\Phi_E^*\cong \Psi_{E^{\vee}\otimes \pi^*\omega_X}[n].
$$

From proposition \ref{comp} the object $K$ of $\db {M\times M}$, which determines $\Phi_E^*\circ\Phi_E$, is 
\begin{equation}\label{K}
K={q_{13}}_*({q_{23}}^*(E^{\vee}\otimes \pi^*\omega_X)\otimes {q_{12}}^*E)[n],
\end{equation}
where the morphisms $q_{13}, q_{23}, q_{12}$ and $\pi$ are taken
from the following diagram

\begin{figure}[th]
\hspace*{5cm}\epsffile{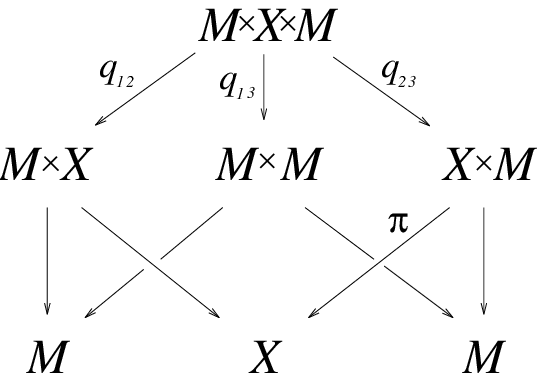}
\end{figure}

We need to prove that $K$ is quasiisomorphic to ${\o\Delta}=\Delta_*{\o M}$, where $\Delta: M\lto M\times M$ is the diagonal embedding, because ${\o\Delta}$ gives the identity functor on $\db{M}$.

Let us consider a commutative diagram 
$$
\begin{array}{ccc}
X&\stackrel{j_{t_1 t_2}}{\longrightarrow}& M\times X\times M\\
\llap{\footnotesize $f$} \downarrow &&\llap{$q_{13}$} \downarrow\\
{\rm Spec} k&\stackrel{i_{t_1 t_2}}{\longrightarrow}& M\times M 
\end{array}
$$
Here $i_{t_1 t_2}$ is the embedding of a geometric point $(t_1 t_2)$ in $M\times M$, and $f : X\to {\rm Spec} k$ the corresponding fibre of $q_{13}$ over this point.

This diagram is useful for computing the fibres of $K$ over points of $M\times M$.

Indeed, by lemma \ref{isof}
$$
i_{t_1 t_2}^*K=i_{t_1 t_2}^*{q_{13}}_*({q_{23}}^*(E^{\vee}\otimes  \pi^*\omega_X)\otimes {q_{12}}^*E)[n]=
$$
\begin{eqnarray}\label{fib}
=f_*j_{t_1 t_2}^*({q_{23}}^*(E^{\vee}\otimes  \pi^*\omega_X)\otimes {q_{12}}^*E)[n]=f_*(j_{t_1 t_2}^*{q_{23}}^*(E^{\vee}\otimes  \pi^*\omega_X)\otimes j_{t_1 t_2}^*{q_{12}}^*E)[n].
\end{eqnarray}
From the commutative diagram

\begin{figure}[h]
\hspace*{5cm}\epsffile{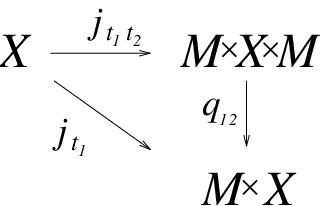}
\end{figure}

where $j_{t_1}$ is the embedding $x\mapsto ( t_1, x)$, and from the definition of $\Phi_E$one obtains: 
\begin{equation}\label{fun1}
j_{t_1 t_2}^*{q_{12}}^*E=j_{t_1}^*E=\Phi_E({\cal O}_{t_1}).
\end{equation}

Analogously,
\begin{equation}\label{fun2}
j_{t_1 t_2}^*{q_{12}}^*(E^{\vee}\otimes  \pi^*\omega_X)=\Phi_E({\cal O}_{t_2})^{\vee}\otimes\omega_X.
\end{equation}
Formulas (\ref{fib}), (\ref{fun1}), (\ref{fun2}) imply isomorphisms:
$$
i_{t_1 t_2}^*K=f_*(\Phi_E({\cal O}_{t_1})\otimes\Phi_E({\cal O}_{t_2})^{\vee}\otimes\omega_X)[n]=
$$
\begin{eqnarray}\label{res}
=f_*({\cal H}om(\Phi_E({\cal O}_{t_2})\:,\;\Phi_E({\cal O}_{t_1}))\otimes\omega_X)[n]={\h \Phi_E({\cal O}_{t_1}), {\Phi_E({\cal O}_{t_2})}}^*.
\end{eqnarray}
The last equality comes from Serre duality on $X$.

Apply proposition \ref{tor} to the diagonal embedding of $M$ in $M\times M$. By formula (\ref{res}) and assumptions of the theorem, the object $K$ satisfies the hypothesis of the proposition. Therefore, $K$ is a pure sheaf with the support at the diagonal $\Delta M$.


The natural transformation $\alpha$ gives rise to a sheaf homomorphism $K\to {\o\Delta}$. It is an epimorphism, because otherwise its image would not generate the stalk of ${\o\Delta}$ at some point $(t,t)$ at the diagonal. But this would imply that $\Phi_E({\cal O}_{\Delta})$ has no endomorphisms ( that is, the trivial object) in contradiction with assumptions of the theorem.

Let $F$ be the kernel of this morphism, i.e. there is an exact sequence of coherent sheaves on $M\times M$ :
\begin{equation}\label{seq}
0\lto F \lto K\lto {\o\Delta}\lto 0
\end{equation}
We have to prove that $F$ is trivial. Considering the pull back of the short exact sequence to any point from $M\times M$ we obtain a long exact sequence showing that the sheaf $F$ satisfies hypothesis of proposition \ref{tor}. It follows from the proposition that the support of $F$ coincides with the diagonal $\Delta M$. It is sufficient to prove that the restriction of $F$ to the diagonal is zero. Let us consider for this the commutative diagram:
$$
\begin{array}{ccc}
M\times X&{\longrightarrow}& M\times X\times M\\
\llap{\footnotesize $p$} \downarrow && \downarrow\\
M&\stackrel{\Delta}{\longrightarrow}& M\times M 
\end{array}
$$
where vertical morphisms are natural projections. Applying lemma \ref{isof} to the object $({q_{23}}^*(E^{\vee}\otimes \pi^*\omega_X)\otimes {q_{12}}^*E)[n]$ from $\db{M\times X\times M}$ and formula (\ref{K}) we obtain a formula for the  derived functors of the restriction--to--diagonal functor for $K$ :
$$
L^i\Delta^*(K)=R^{n-i} p_*(E\otimes E^{\vee}\otimes \pi^*\omega_X).
$$
Therefore, by the relative version of Serre duality and hypothesis of the theorem
$$
\Delta^*K={\o \Delta},
$$
$$
L^1 \Delta^*(K)=R^1p_*(E\otimes E^{\vee})^{\vee}.
$$

Unfortunately, it is not sufficient to know that the restriction of $K$ to the diagonal is $ {\o \Delta}$, because $K$ might not be the push forward along $\Delta$ of a sheaf on $M$ (being, `situated' on some infinitesimal neighborhood of  $\Delta(M)$).
 
Furthermore, $L^1 \Delta^*({\o \Delta})=\Omega^1_M$, this means that the long exact sequence, obtained from (\ref{seq}) by tensoring with ${\o \Delta}$  looks as follows:
$$
\cdots\lto R^1p_*(E\otimes E^{\vee})^{\vee}\stackrel{\beta}{\lto} \Omega^1_M\lto L^0\Delta^* F\lto{\o \Delta}\stackrel{\sim}{\lto}{\o \Delta}\lto 0.
$$ 

Since the support of $F$ coincides with $\Delta(M)$, so does the support of $L^0\Delta^* F$, i.e. $L^0\Delta^* F$ is not a torsion sheaf.

Therefore, if $L^0\Delta^* F$ is not zero, then this exact sequence shows that
$$
\beta^* : T_M\lto R^1p_*(E\otimes E^{\vee})
$$
has a non-trivial kernel.
\par
\bigskip

{\sc Remark.}
 If one consider $\Phi_E({\cal O}_t)$ as a system of objects from $\db{X}$ parametrized by $M$, then the restriction of $\beta^*$ to any point $t$ from $M$ is , actually, the homomorphism from the deformation theory
$$
T_{M,t}\lto {\H 1, X, {\Phi_E({\cal O}_t)}, {\Phi_E({\cal O}_t)}}.
$$

Therefore, a vector field from the kernel of $\beta^*$ gives a direction, which the objects do not change along with. This is in contradiction with the orthogonality assumptions of the theorem. Unfortunately, integrating such an algebraic vector field one might obtain non-algebraic curves. 
\bigskip

For this reason our further strategy is going to find only a formal one--parameter deformation at one point $t_0$ in $M$, along which $E$ has a formal connexion (analog of trivialization), and then to bring this in contradiction with  the property of $\Phi_E({\cal O}_{t_0})$ to having the support in point $t_0$, which is a consequence of the orthogonality condition.  

Consider a point $t_0$ in $M$, $U$ an open neighborhood of 
$t_0$ and a non-zero at $t_0$ local section $\xi\in H^0(U , T_M\Bigl|_U)$, which belongs to the kernel of $\beta^*$. This vector field $\xi$ defines a formal 1--dimensional subscheme 
$\mit\Gamma$ of the formal neighborhood $\hat U$ of $t_0$
in $M$. The defining ideal of the subscheme consists of the  function on ${\hat U}$ having trivial all iterated derivatives along $\xi$ at point $t_0$ ( the zero derivative being the value of a function at $t_0$ ) :
$$
I=\Bigl\{ f\in H^0(\hat U, {\cal O})\; \Bigl| \;\xi^k(f)|_{t_o}=0, \:\mbox{for any }\: k\ge 0 \Bigl\}.
$$ 

It follows that the restriction of $\beta^*$ to the tangent bundle $T_{\mit\Gamma , t_0}$ of $\mit\Gamma$ at $t_0$ is trivial.

Denote $E_{\mit\Gamma}$ the restriction of $E$ to $\mit\Gamma\times X$. One has

\begin{eqnarray}\label{R}
{\H 1, {\mit\Gamma\times X}, {p^*T_{\mit\Gamma} }, {E_{\mit\Gamma}^{\vee}\otimes E_{\mit\Gamma}}}\cong 
{\H 1, {\mit\Gamma}, {T_{\mit\Gamma}} , {p_*(E_{\mit\Gamma}^{\vee}\otimes E_{\mit\Gamma})}}\cong\nonumber\\
{\H 0 , {\mit\Gamma}, {T_{\mit\Gamma}}, {R^1p_*(E_{\mit\Gamma}^{\vee}\otimes E_{\mit\Gamma})}},
\end{eqnarray}
since $T_{\mit\Gamma}$ is free ( of rank 1 ) on ${\mit\Gamma}$.

Let us consider the first infinitesimal neighborhood $\Delta^1_{\mit\Gamma}$ of the diagonal $\Delta_{\mit\Gamma} : {\mit\Gamma}\times X\lto {\mit\Gamma}\times {\mit\Gamma}\times X$. Pulling $E_{\mit\Gamma}$ back to ${\mit\Gamma}\times {\mit\Gamma}\times X$ 
 along the first coordinate, then restricting to $\Delta^1_{\mit\Gamma}$ and then pushing forward along the second coordinate, one obtains the object $J^1(E_{\mit\Gamma})\in \db{{\mit\Gamma}\times X}$ of `first jets' of $E_{\mit\Gamma}$. It is included in an exact triangle :
\begin{equation}\label{tr}
E_{\mit\Gamma}\otimes p^*\Omega^1_{\mit\Gamma}\lto J^1(E_{\mit\Gamma})\lto E_{\mit\Gamma}\stackrel{at_E}{\lto}
E_{\mit\Gamma}\otimes p^*\Omega^1_{\mit\Gamma}[1].
\end{equation}
Here $at_E$ is the so called Atiyah class of $E$. It can be considered as an element of ${\H 1, {\mit\Gamma\times X}, {p^*T_{\mit\Gamma} }, {E_{\mit\Gamma}^{\vee}\otimes E_{\mit\Gamma}}}$. Under identification from (\ref{R}) $at_E$ comes into the restriction of $\beta^*$ to $\mit\Gamma$, which is the trivial element of ${\H 0 , {\mit\Gamma}, {T_{\mit\Gamma}}, {R^1p_*(E_{\mit\Gamma}^{\vee}\otimes E_{\mit\Gamma})}}$ by the choice of $\mit\Gamma$.

We consider $E_{\mit\Gamma}$ as an element  of the derived category of quasicoherent sheaves on X. It is naturally endowed with an additional homomorphism $A\to {\rm End}_X E_{\mit\Gamma}$, where $A$ is an algebra of functions on ${\mit\Gamma}$ ( isomorphic to $k[[t]]$). Such a homomorphism          we call by {\it A}--module structure on $E_{\mit\Gamma}$. An {\it A}--module structure on $E_{\mit\Gamma}$ induces {\it A}--module structures on $E_{\mit\Gamma}\otimes p^*\Omega^1_{\mit\Gamma}$ and $J^1(E_{\mit\Gamma})$, so that morphisms from (\ref{tr}) are compatible with them.

Like for usual vector bundles there exists a natural homomorphism in $D^b_{Qcoh}(X)$ :
$$
 E_{\mit\Gamma}\stackrel{\mu}{\lto}J^1(E_{\mit\Gamma}) , 
$$
which is a differential operator of the first order with respect to the {\it A}--module structures. Triviality of $at_E$ implies existence of a morphism 
$$
J^1(E_{\mit\Gamma})\stackrel{\nu}{\lto}E_{\mit\Gamma}\otimes p^*\Omega^1_{\mit\Gamma},
$$
which is a section of the first morphism from (\ref{tr}). The composition $\nabla=\nu\circ\mu$ defines a morphism of quasicoherent sheaves on $\mit\Gamma\times X$
$$
\nabla : E_{\mit\Gamma}\lto E_{\mit\Gamma}\otimes p^*\Omega^1_{\mit\Gamma},
$$
which is a connexion on $E_{\mit\Gamma}$ along the fibres of the projection $p_{\mit\Gamma} :{\mit\Gamma}\times X \lto {\mit\Gamma}$ in the sense that if $t\in A$ is a function on our formal scheme $\mit\Gamma$ , then the following equality for morphisms from $E_{\mit\Gamma}$ to $E_{\mit\Gamma}\otimes p^*\Omega^1_{\mit\Gamma}$ is valid:
\begin{equation}\label{con}
\nabla{\circ}t - t{\circ}\nabla=dt, 
\end{equation}
here $t$ is identified with the corresponding morphism from $E_{\mit\Gamma}$ to $E_{\mit\Gamma}$ and $dt$ denotes the operator of tensor multiplication by $dt$.

Since ${\mit\Gamma}$ is a one-dimensional subscheme, $\Omega^1_{\mit\Gamma}$ is a one-dimensional free {\it A}--module. Hence, for the reason of simplicity we can 
identify $E_{\mit\Gamma}$ with $E_{\mit\Gamma}\otimes p^*\Omega^1_{\mit\Gamma}$ by means of tensoring with $dt$, where $t\in A$ is a formal parameter on the scheme ${\mit\Gamma}$. 
Then, formula (\ref{con}) gives the coordinate-impulse relation between $\nabla$ and $t$:
\begin{equation}\label{ci}
[\nabla, t ]={\bf 1}.
\end{equation}

\th{Lemma}\label{inv}
Under the above identification of $E_{\mit\Gamma}$ with $E_{\mit\Gamma}\otimes p^*\Omega^1_{\mit\Gamma}$ the morphism $\nabla{\circ}t$ is invertible in ${\rm End}_X E_{\mit\Gamma}$
\eth
\pr
From (\ref{ci}) one has :
$$
[\nabla, t ^k]=k t^{k-1}.
$$
This gives a formula for the inverse to $\nabla\circ t $ :
$$
(\nabla{\circ}t)^{-1}=\sum_{k=0}^{\infty}\frac{(-1)^k t^k{\circ}\nabla^k}{(k+1)!}
$$
This formal series correctly defines an endomorphism of $E_{\mit\Gamma}$, because by definition $E_{\mit\Gamma}$ is the limit of a system of objects $E_{\mit\Gamma_n}$ from 
${\mit\Gamma}_n\times X$ , where ${\mit\Gamma}_n$ is the n--th infinitesimal neighborhood of $t_0$ in ${\mit\Gamma}$. For every $n$ the formula gives a finite expansion for an  endomorphism of $E_{\mit\Gamma_n}$, thus , in the limit, it does an endomorphism of $E_{\mit\Gamma}$ . This proves the lemma.
\bigskip

Let $E_0$ be the first member of the exact triangle :
$$
E_0\stackrel{\rho}{\lto}E_{\mit\Gamma}\stackrel{\nabla}{\lto}
E_{\mit\Gamma}.
$$
It is an object of the derived category of quasicoherent ${\o X}$--modules (`horisontal sections of $E$').

\th{Proposition}\label{flat}
i) The composition $\lambda$ of $\rho\times id_A$ with the multiplication morphism  $E_{\mit\Gamma}\times A\lto E_{\mit\Gamma}$:
$$
\lambda : E_0\times A \stackrel{\rho\times id_A}{\lto}E_{\mit\Gamma}\times A\lto E_{\mit\Gamma}
$$
is an isomorphism, in other words it yields trivialization of $E_{\mit\Gamma}$.

ii) $E_0$ is quasiisomorphic to a complex of coherent sheaves on $X$.
\eth
\pr
Let us consider the cone $C$ of $\lambda$ :
$$
E_0\times A \stackrel{\lambda}{\lto}E_{\mit\Gamma}\lto C.
$$

Restricting this exact triangle to the fibre $X_0$ of $p_{\mit\Gamma}$ over the closed point of ${\mit\Gamma}$  (which is, of course, naturally identified with $X$), one obtains an exact triangle :
$$
E_0\stackrel{\lambda_0}{\lto}E_{\mit\Gamma}\Bigl|_{X_0}\lto C\Bigl|_{X_0}.
$$
Vanishing of $C\Bigl|_{X_0}$ implies vanishing of $C$, hence, for proving i) we need to show that the left morphism $\lambda_0$ of this triangle is isomorphism.

Multiplication by $t$ gives an exact triangle of sheaves on $
\mit\Gamma\times X$ :
$$
0\lto{\o {\mit\Gamma\times X}}\stackrel{t}{\lto}{\o {\mit\Gamma\times X}}\lto{\o {X_0}}\lto 0
$$
It lifts to an exact triangle:
$$
E_{\mit\Gamma}\stackrel{t}{\lto}E_{\mit\Gamma}\lto E_{\mit\Gamma}\Bigl|_{X_0}\lto E_{\mit\Gamma}[1].
$$
Consider an octahedral diagram of exact triangles \cite{BBD}:

\begin{figure}[th]
\hspace*{5cm}\epsffile{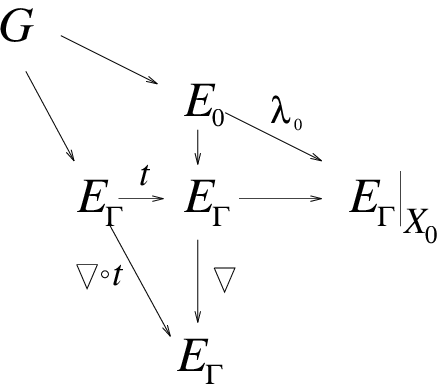}
\end{figure}

By lemma \ref{inv} $\nabla{\circ}t$ is an isomorphism. Hence, $G$ is zero object and $\lambda_0$ is an isomorphism. Since $ E\Bigl|_{X_0}$ is the restriction of complex of coherent sheaves to $X_0$, $E_0$ is coherent over $X_0$. Since $C\Bigl|_{X_0}=G[1]$ is zero, so is $C$. Therefore, $\lambda$ is isomorphism.

In order to finish the proof of theorem \ref{mai} let us look at the image $L=\Phi^*_E \circ \Phi_E({\o {t_0}})$ of ${\o {t_0}}$ under the functor
$$
\Phi^*_E\circ \Phi_E: \db{M}\lto\db{M}.
$$
Recall that by lemma \ref{adj}
$$ 
\Phi^*_E=\Psi_{E^{\vee}\otimes \pi^*\omega_X}[n].
$$

The trivialization of $E$ along ${\mit\Gamma}$ from proposition \ref{flat} gives us a similar trivialization of $E^{\vee}\otimes \pi^*\omega_X$. By the definition of $\Psi$ this implies trivialization along ${\mit\Gamma}$ of any object from the image of $\Phi^*_E$. 
Since we know that $\Phi^*_E\circ \Phi_E$ is determined by sheaf $K$, having the  diagonal $\Delta(M)$ as its support, the image $L$ of a skyscraper sheaf ${\o {t_0}}$ is a non--zero object from $\db{M}$ having $t_0$ as the support.

This means that $L$ annihilates by some power $I^k$ of the  maximal ideal $I\subset A$. Such an object has a trivialization only if it is zero. This finishes the proof of theorem \ref{mai}.
\bigskip


The simplest example of a full and faithful functor $\db{M}\lto\db{X}$ rises in the case when $M$ is a point. In this situation we have the only {\sf object} $E\in \db{X}$, which is an {\sf exceptional} one:
$$
\begin{array}{ll}
{\H 0, X, E, E}=k, &\\
{\H i, X, E, E}=0, \quad\mbox{for }\: i\ne 0
\end{array}
$$

It gives a functor from the derived category of vector spaces over $k$ to $\db{X}$.


Mukai in \cite{Mu} and \cite{Mu2} considered two important examples of fully faithful functors between geometric categories.

First one is the so called Fourier--Mukai transform. It gives equivalence
$$
\db{A}\lto\db{\hat A}
$$
for any abelian variety $A$ and its dual $\hat A$.

We briefly recall his construction.

Let $A$ be an abelian variety of dimension $g$, $\hat A$ its dual abelian variety and ${\cal P}$ the normalized Poincare bundle on $A\times\hat A$. As $\hat A$ is a moduli space of invertible sheaves on $A$, ${\cal P}$ is a linear vector bundle, and normalization means that both ${\cal P}\Bigl|_{A\times\hat 0}$ and ${\cal P}\Bigl|_{0\times\hat A}$ are trivial.
\th{Theorem}{\rm (\cite{Mu}).}\label{ab}
The functors $\Phi_{\cal P}:\db{A}\to \db{\hat A}$ and $\Psi_{\cal P}:\db{\hat A}\to \db{A}$ are equivalences of triangulated categories and
$$
\Psi_{\cal P}\circ\Phi_{\cal P}\cong (-1_A)^*[g],
$$
$$
\Phi_{\cal P}\circ\Psi_{\cal P}\cong (-1_{\hat A})^*[-g],
$$
here $(-1_A)^*$ is the auto-equivalence of $\db{A}$ induced by the automorphism of multiplication by $-1$ on $A$.
\eth
\pr (see \cite{Mu}).

In the case of a principally polarized abelian variety $(A, L)$, where $L$ is a polarization, the dual $\hat A$ is identified with $A$. Then $\Phi_{\cal P}$ can be regarded as an auto-equivalence of $\db{A}$. $\Phi_{\cal P}$ in couple
with the functor of tensoring by $L$ generates the action of 
the Artin braid group $B_3$ on three strands.

The other example of Mukai is a K3--surface $S$, while $M$ is a moduli space of stable vector bundles.

Specifically, for a smooth  K3--surface $S$ one consider the Mukai lattice ${\cal M}(S)$, which is the image of the Chern homomorphism $K_0(S)\lto H^*(S, {\bf C})$ from the Grothendieck
group $K_0(S)$ to full cohomology group $H^*(S, {\bf C})$. There is the Euler bilinear form on ${\cal M}(S)$, which for vectors $v$ and $v'$ presented by some sheaves ${\cal F}$ and
 ${\cal F'}$ is defined by the formula:
$$
\chi(v, v')= \sum (-1)^i {\rm dimExt}^i({\cal F, F'}).
$$
Since the canonical class is trivial, by Serre duality this form is symmetric.

Let $v$ be an isotropic indivisible by integer vector with respect to $\chi$. The coarse moduli space of stable bundles on $S$, corresponding to $v$, is again a smooth K3--surface $S'$. There is a rational correspondence between $S$ and $S'$. If $S'$ is a {\it fine}  moduli space, then we have the universal vector bundle $E$ on $S\times S'$.
\th{Theorem}{\rm \cite{Mu2}}. 
Functor $\Phi_{E}:\db{S}\lto \db{S'}$ is an equivalence of triangulated categories.
\eth

In the both examples of equivalences the canonical class of varieties (either of abelian one or of a K3--surface) is trivial. In chapter 3 we construct another example of equivalence between geometric categories using flops. The centre of such transformation is in a sense trivial with respect to the canonical class. An explanation for this phenomenon
is given in chapter 4.

\section{Intersection of two even dimensional quadrics.}

In this chapter we show how theorem 1.1 helps to construct a semiorthogonal decomposition of the derived category of coherent sheaves on the intersection of two even dimensional quadrics, with one summand
being the derived category on a hyperelliptic curve and with the others being generated by single exceptional objects.

This result can be considered as a categorical explanation 
for the description, due to Desale and Ramanan, of moduli spaces of rank 2 vector bundles on a hyperelliptic curve as a base of a family of projective subspaces belonging to the intersection of two even dimensional quadrics \cite{DR}. Our construction gives analogous description for any moduli spaces of bundles on the curve by means of families of complexes of coherent sheaves on the intersection locus.

We first recall some definitions and facts concerning exceptional sequences, admissible subcategories, Serre functors and semiorthogonal decompositions \cite{Bon}, \cite{BK}.

Let ${\cal B}$ be a full subcategory of an additive category. The {\sf right orthogonal} to ${\cal B}$ is the full subcategory ${\cal B}^{\perp}\subset {\cal A}$ consisting of the objects $C$ such that ${\h B, C}=0$
for all $B\in{\cal B}$. The {\sf left orthogonal} ${}^{\perp}{\cal B}$ is defined analogously. If ${\cal B}$ is a triangulated subcategory of a triangulated category ${\cal A}$, then ${}^{\perp}{\cal B}$ and ${\cal B}^{\perp}$ are also triangulated subcategories.
\th{Definition}\label{adm}
Let ${\cal B}$ be a strictly full triangulated subcategory of  a triangulated category ${\cal A}$. We say that ${\cal B}$ is {\sf right admissible} (resp., {\sf left admissible}) if for each $X\in{\cal A}$ there is an exact triangle $B\to X\to C$, where $B\in{\cal B}$ and $C\in{\cal B}^{\perp}$ (resp., $D\to X\to B$, where $D\in{}^{\perp}{\cal B}$ and $B\in{\cal B}$). A subcategory is called {\sf admissible} if it is left and right  admissible.
\eth

\th{Definition}\label{exc}
An {\sf exceptional object} in a derived category ${\cal A}$ is an object $E$ satisfying the conditions ${\rm Hom}^i(E\:,\;E)=0$ when $i\ne0$ and ${\rm Hom}(E\:,\;E)=k$.
\eth
\th{Definition}\label{exs}
{\sf A full exceptional sequence} in ${\cal A}$ is a
sequence of exceptional objects $(E_0,..., E_n)$, satisfying the semiorthogonal condition ${\rm Hom}^.(E_i\:,\;E_j)=0$ when $i>j$, and generating the category ${\cal A}$.
\eth

The concept of an exceptional sequence is a special case of the concept of a semiorthogonal sequence of subcategories:

\th{Definition}\label{sd}
A sequence of admissible subcategories $({\cal B}_0,..., {\cal B}_n)$ in a derived category ${\cal A}$ is said to be {\sf semiorthogonal} if the condition ${\cal B}_j\subset {\cal B}^{\perp}_i$ holds when $j<i$ for any $0\le i\le n$. In addition, a semiorthogonal sequence is said to be {\sf full} if it generates the category ${\cal A}$. In this case we call such a sequence {\sf semiorthogonal decomposition} of the category ${\cal A}$ and denote this as follows:
$$
{\cal A}=\Bigl\langle{\cal B}_0,....,{\cal B}_n\Bigl\rangle.
$$
\eth

\th{Definition}\label{SF}
Let ${\cal A}$ be a triangulated $k$--linear category with finite--dimensional ${\rm Hom's}$. A covariant additive functor $F: {\cal A}\to{\cal A}$ that commutes with translations is called a {\sf Serre functor} if it is a category equivalence and there are given bi--functorial isomorphisms
$$
\varphi_{E,G}: {\rm Hom}_{\cal A}(E\:, \;G)\stackrel{\sim}{\lto}{\rm Hom}_{\cal A}(G\:, \;F(E))^*
$$
for $E,G\in{\cal A}$, with the following property: the composite
$$
(\varphi_{F(E),F(G)}^{-1})^*{\scriptsize{\circ}}\varphi_{E,G}: {\rm Hom}_{\cal A}(E\:, \;G)\lto{\rm Hom}_{\cal A}(G\:, \;F(E))^*\lto{\rm Hom}_{\cal A}(F(E)\:, \;F(G))
$$
coincides with the isomorphism induced by $F$.
\eth

\th{Theorem}\label{USF}{\rm \cite{BK}}
i) Any Serre functor is exact,

ii) Any two Serre functors are connected by a canonical
functorial isomorphism.
\eth
Let $X$ be a smooth algebraic variety, $n={\rm dim}X$, ${\cal A}=\db{X}$ the derived category of coherent sheaves on $X$, and $\omega_X$ the canonical sheaf. Then the functor $(\cdot)\otimes\omega_X[n]$ is a Serre functor on ${\cal A}$, in view of the Serre--Grothendieck duality:
$$
{\rm Ext}^i(F\:,\;G)={\rm Ext}^{n-i}(G\:,\;F\otimes\omega_X)^* 
$$

Let us fix notations. For vector spaces $U$ and $V$ of dimension 2 and $n=2k$, respectively, we consider a linear embedding:
$$U\stackrel{\varphi}{\longrightarrow} S^2V^*.$$

By projectivization $\varphi$ defines a pencil of projective quadrics in ${\bf P}^{n-1}={\bf P}(V)$, parametrized by ${\bf P}^1={\bf P}(U)$.

Denote by $X$ the intersection locus of these quadrics. Let 
$\{q_i\}_{i=1,\ldots,n}\subset {\bf P}^1$ are the points, corresponding to the degenerate quadrics. We assume that all $q_i$ are mutually distinct. This implies that $X$ is a smooth variety and quadrics corresponding to $q_i$ have simple degeneration. Consider a double covering 
$C\stackrel{p}{\longrightarrow} {\bf P}^1$
with ramification in all points $\{q_i\}$. Then $C$ is a hyperelliptic curve.

In order to construct a fully faithful functor
$D_{coh}^b(C)\longrightarrow D_{coh}^b(X)$ we find a vector bundle $S$ on $C\times X$ and then use theorem 1.1 for the functor $\Phi_S$ (see ch.1). To outline the idea of constructing the bundle $S$, let us recall that for non-degenerate even dimensional quadric there exist two spinor bundles (c.g.~\cite{KAPR}). Restricting these two bundles  to $X$  and varying our quadric in the pencil we obtain that $C$ is the fine moduli space of spinor bundles. 
Unfortunately, the fine moduli space exists only for a pencil of even dimensional quadrics. For the case of more than two quadrics of arbitrary dimension there appear some global obstructions for gluing together spinor bundles and local problems for extending to points, corresponding to degenerate quadrics.

A generalization to the case of more then two quadrics of arbitrary  dimension will be given in a forthcoming paper.

Let $Y$ (relative grassmanian of maximal isotropic subspaces)
be a subvariety in ${\bf P}(U)\times {\rm G}(k,V)$ consisting of the pairs $(q,L)$ such that $L$ is isotropic with respect to the quadric corresponding to $q$ (which we denote by the same letter $q$):
$$Y:=\Bigl\{ (q,L)\in{\bf P}^1\times {\rm G}(k,V)|\ q(L)=0
\Bigr\}.$$
The image of $Y$ under the natural map into Albanese variety is isomorphic to $C$. Thus, we have a natural projection 
$\varphi:Y\longrightarrow C$,
which is a smooth projective morphism. Its fibre over a point $c\in C$ is one of two connected components of the maximal isotropic grassmanian, corresponding to the quadric $p(c)$. Of course, the composition 
$p{\scriptstyle\circ}\varphi$
coincides with the natural projection to ${\bf P}^1$, in other words $p$ and $\varphi$ give the Stein factorization of the projection.

Now consider linear subspaces of dimension $k-1$, which belongs to $X$. It is well known that the variety of all such subspaces is isomorphic to Jacobian $J(C)$ of the curve $C$
 \cite{R1}. 
 We choose one of them $M$. It gives a section of $\varphi$. Indeed, if one consider a subvariety $C_M\subset Y$ of pairs $(q,L)\in Y$ such that $L$ contains $M$:
$$C_M:=\Bigl\{ (q,L)\in Y|\  L\supset M
\Bigr\},$$
then $\varphi$ biregulary projects $C_M$ to $C$, because for any non-degenerate (resp., degenerate) quadric from our pencil there exist two (resp., one) containing $M$ maximal isotropic subspaces, which lie in the different components of the grassmanian, corresponding to this quadric.

Now consider the subvariety $D\subset Y$ of pairs $(q,L)$ such that $L$ has a non-trivial intersection with $M$:
$$D:=\Bigl\{ (q,L)\in Y|\  L\cap M\not=0
\Bigr\}.$$
Then $D$ is a divisor in $Y$.
Denote by ${\cal L}={\cal O}(D)$ the corresponding linear bundle on $Y$. 

Now consider the variety $F$ (of partial isotropic flags) consisting of triples $(q,l,L)\in {\bf P}^1\times {\bf P}(V)\times{\rm G}(k,V)$ such that $l\subset L$ and $q(L)=0$:
$$F:=\Bigl\{ (q,l,L)\in {\bf P}^1\times {\bf P}(V)\times{\rm G}(k,V)|\  l\subset L,\ q(L)=0\Bigr\}.$$
Since $l\subset L$, $l$ is a point of the quadric $q$. In other words, the projections of $F$ to the components
${\bf P}^1\times{\bf P}(V)$ and ${\bf P}^1\times{\rm G}(k,V)$ of the product give the couple of maps: 
$\mu:F\longrightarrow Y,\ \lambda:F\longrightarrow Q,$ where $Q$ is the relative quadric, i.e., the variety of pairs $(q,l)\in{\bf P}^1\times{\bf P}(V)$, such that $q(l)=0$:
$$Q:=\Bigl\{ (q,l)\in {\bf P}^1\times {\bf P}(V)|\   q(l)=0\Bigr\}.$$
The projection to the first component gives a map: $Q\longrightarrow{\bf P}^1$.

Let $Q_C=Q\times_{{\bf P}^1}C$ be the product of $Q$ over ${\bf P}^1$. Since $\varphi{\scriptstyle\circ}\mu:F\longrightarrow C$ and $\lambda:F\longrightarrow Q$ are maps to the component of the product, by the universal property we have a map $\nu:F\longrightarrow Q_C$.

Since $X$ belongs to all quadrics of the pencil we have the natural embedding $X\times{\bf P}^1\hookrightarrow Q$, which lifts up to an embedding $\varepsilon:X\times C\hookrightarrow Q_C$. All these varieties and maps are depicted in the following diagram:

\centerline{\epsffile{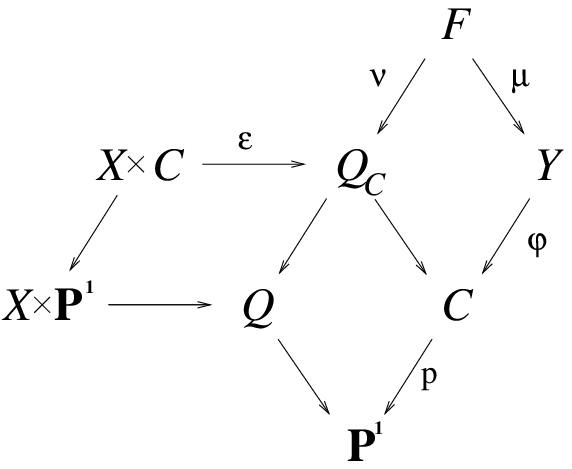}}

Now define $S=\varepsilon^*\nu_*\mu^*{\cal L}.$ Actually, $S$ is a vector bundle on $X\times C$. 

Let us fixed a point $c\in C$. If $q=p(c)$ is a smooth quadric, then the fibre of $S$ over $X\times c\simeq X$ is one of two spinor bundles on $q$, restricted to $X$. If 
$q=p(c)$ is degenerate, then it is a cone over a quadric of 
the same  dimension as $X$. Then the fibre of $S$ over $X\times c\simeq X$ is the restriction to $X$ of the pull back of the spinor bundle on this even dimensional quadric to the cone. Since $X$ does not meet the singular point of the cone, this restriction is also a vector bundle on $X$.

Let us recall the structure of the derived category for a smooth projective quadric due to M.~Kapranov.

There exist two for an even dimensional (resp., one for an odd dimensional ) quadric $q$ spinor bundles $S_q$ and $\widetilde S_q$ (resp., $S_q$). The exceptional sequence
\begin{equation}\label{qcol}
\Bigl({\cal O}(-d+1),{\cal O}(-d+2),\ldots,{\cal O},\widetilde
S_q,S_q\Bigr)\qquad{\rm for}\  d\ {\rm even},
\end{equation}
\begin{equation}\label{qcol2}
\Bigl({\cal O}(-d+1),{\cal O}(-d+2),\ldots,{\cal O},S_q\Bigr)\qquad\qquad{\rm for}\  d\ {\rm odd}
\end{equation}
is a full strong exceptional sequence on $q$, here $d=\dim q$
(see \cite{KAPR}).  For $d$ even, $S_q$ and $\widetilde S_q$ are mutually homologically orthogonal.

\th{Theorem}\label{iq} The functor $\Phi_S:D_{coh}^b(C)\longrightarrow D_{coh}^b(X)$ is full and faithful.
\eth
\pr
We shall use facts about spinor bundles on smooth quadrics.

We have to verify conditions of theorem 1.1 for fibres of $S$ over $X\times c$.

First, let us check the orthogonality conditions. Suppose that $c_1$ and $c_2$ are points of $C$, such that $p(c_1)\not=p(c_2)$, and let $S_1$ and $S_2$ are spinor bundles over corresponding quadrics. There are short exact sequences of sheaves on the projective space ${\bf P}(V)$:
$$
0\longrightarrow V\otimes{\cal O}(-2)\longrightarrow 
V\otimes{\cal O}(-1)\longrightarrow S^{\vee}_1\longrightarrow 0,
$$
$$
0\longrightarrow W\otimes{\cal O}(-1)\longrightarrow W\otimes {\cal O}\longrightarrow S_2 \longrightarrow 0.
$$
Here we identify bundles on quadrics with corresponding coherent sheaves on ${\bf P}(V)$. If any of these quadrics is degenerate, then the same sequences holds beyond singular points  of the quadric, being sufficient for what follows. Consider these sequences as resolutions for $S^{\vee}_1$, and $S_2$ and use them for computation $S_1^{\vee}\otimes S_2$. Since the quadrics intersect transversally, there are no torsion groups:
$$
{\rm Tor}^i_{{\bf P}(V)}(S_1^{\vee},S_2)=0,\qquad{\rm for}\ i>0.
$$
Therefore, we obtain a resolution for $S_1^{\vee}\otimes S_2$ of the following kind:
$$
0\longrightarrow C\otimes {\cal O}(-3)\longrightarrow 
 B\otimes {\cal O}(-2)\longrightarrow 
 A\otimes {\cal O}(-1)\longrightarrow S_1^{\vee}\otimes S_2\longrightarrow 0.
$$
Computing cohomologies of $S_1^{\vee}\otimes S_2$ by means of this resolution we obtain the orthogonality conditions for the case $p(c_1)\not=p(c_2):$
$$
{\rm Ext}^i_X(S_1,S_2)={\rm H}^i({\bf P}(V),S_1^{\vee}\otimes S_2)=0.
$$

Now suppose that $p(c_1)=p(c_2)$. Then we have to verify orthogonality, while restricted to $X$, between two spinor bundles $S_q$ and $\widetilde S_q$ on a single non-degenerate quadric $q$.

Consider the tensor product $S_q^{\vee}\otimes\widetilde S_q$ over $q$. Since $X$, as a divisor in $q$, is equivalent to double hyperplane section, we have an exact sequence of sheaves on $q$:
$$
0\longrightarrow S_q^{\vee}\otimes\widetilde S_q(-2)\longrightarrow
S_q^{\vee}\otimes\widetilde S_q \longrightarrow
S_q^{\vee}\otimes\widetilde S_q \left|_X\right. \longrightarrow 0.
$$

Computing cohomology we easily find that 
$$
{\rm Ext}^i_X(S_q,\widetilde S_q)={\rm H}^i(q,S_q^{\vee}\otimes\widetilde S_q \left|_X\right.)=0,\qquad{\rm for\ any}\ i.
$$
Indeed,
$$
{\rm H}^i(S_q^{\vee}\otimes\widetilde S_q )={\rm Ext}^i_X(S_q,\widetilde S_q)=0,\qquad{\rm for\ any}\ i.
$$
Then, mutating $\widetilde S_q$ two times to the left in   sequence (\ref{qcol}) one can obtain a new exceptional sequence:
$$\Bigl({\cal O}(-d+1),\ldots,{\cal O}(-2),\widetilde
S_q(-2),{\cal O}(-1),{\cal O},S_q\Bigr).$$
It yields:
$$
{\rm H}^i(S^{\vee}_q\otimes \widetilde S_q(-2))={\rm Ext}^i(S_q,\widetilde S_q(-2))=0.
$$

Now let us verify condition ii) of theorem 1.1. Suppose $q=p(c)$ is a non-degenerate quadric and $\Phi_S({\cal O}_C)=\widetilde S_q$. We can calculate 
${\rm Ext}^i_X(\Phi_S({\cal O}_C),\Phi_S({\cal O}_C)) $ using the following exact sequence on the quadric:
$$
0 \longrightarrow \widetilde S^{\vee}_q\otimes \widetilde S_q(-2)
 \longrightarrow \widetilde S^{\vee}_q\otimes \widetilde S_q \longrightarrow \widetilde S^{\vee}_q\otimes \widetilde S_q\left |_X\right.\longrightarrow 0.
$$
Since $\widetilde S_q$ is an exceptional object in $D^b_{coh}(q)$ and $\widetilde S_q(-2)$ is a double left mutation in the collection (\ref{qcol}), we have

$$
{\rm H}^0_q(\widetilde S^{\vee}_q\otimes \widetilde S_q)=k,\quad
{\rm H}^i_q(\widetilde S^{\vee}_q\otimes \widetilde S_q)=0\ i\not=0;
$$
$$
{\rm H}^2_q(\widetilde S^{\vee}_q\otimes \widetilde S_q(-2))=k,\quad
{\rm H}^i_q(\widetilde S^{\vee}_q\otimes \widetilde S_q(-2))=0\ i\not=2.
$$
Then the short sequence gives:
$$
{\rm Hom}_X( \widetilde S_q, \widetilde S_q)=k,\ \ 
{\rm Ext}^1_X( \widetilde S_q, \widetilde S_q)=k,\ \ 
{\rm Ext}^i_X( \widetilde S_q, \widetilde S_q)=0,\  i>1.
$$
Similarly for $\Phi_S({\cal O}_C)=S_q$.

Now suppose that $q=p(c)$ is a degenerate quadric. Then the  projection from the centre of the cone gives a double covering $\pi:X\longrightarrow q'$ from $X$ to a quadric of dimension $\ d-1$. Since $\Phi_S({\cal O}_C)=\pi^*S_{q'}$ is the pull back of the spinor bundle $S_{q'}$ on this quadric along $\pi$, 
\begin{equation}\label{sq}
{\rm Ext}^1(\Phi_S({\cal O}_C),\Phi_S({\cal O}_C))=
{\rm H}^i(X,\pi^*(S^{\vee}_{q'}\otimes S_{q'}))=
{\rm H}^i(q',\pi_*\pi^*(S^{\vee}_{q'}\otimes S_{q'})).
\end{equation}
By projection formula we have:
$$
\pi_*\pi^*(S^{\vee}_{q'}\otimes S_{q'})=\pi_*{\cal O}_X\otimes S^{\vee}_{q'}\otimes S_{q'}=\Bigl[ {\cal O}_{q'}\oplus {\cal O}_{q'}(-1)\Bigr]\otimes S^{\vee}_{q'}\otimes S_{q'}=
$$
\begin{equation}\label{pp}
=S^{\vee}_{q'}\otimes S_{q'}\oplus S^{\vee}_{q'}\otimes S_{q'}(-1).
\end{equation}
Since $S_{q'}$ is exceptional on $q'$ and $S_{q'}(-1)$ is the left mutation of $S_{q'}$ in sequence (\ref{qcol2}), it follows that
$$
{\rm H}^0(S^{\vee}_{q'}\otimes S_{q'})=k,\quad 
{\rm H}^i(S^{\vee}_{q'}\otimes S_{q'})=0 \ \ {\rm for}\ i\not=0;
$$
$$
{\rm H}^1(S^{\vee}_{q'}\otimes S_{q'}(-1))=k,\quad 
{\rm H}^i(S^{\vee}_{q'}\otimes S_{q'})=0 \ \ {\rm for}\ i\not=1.
$$
Combining this with (\ref{sq}) and (\ref{pp}), we obtain:
$$
{\rm Hom}\Bigl(\Phi_S({\cal O}_C),\Phi_S({\cal O}_C)\Bigr)=k,\qquad {\rm Ext}^1\Bigl(\Phi_S({\cal O}_C),\Phi_S({\cal O}_C)\Bigr)=k,
$$
$$
{\rm Ext}^i\Bigl(\Phi_S({\cal O}_C),\Phi_S({\cal O}_C)\Bigr)=0,\quad{\rm for}\ i>1.
$$
This concludes the proof of the theorem.

Let us recall that we consider the intersection $X$ of two quadrics of dimension $d=n-2.$

\th{Proposition}\label{orthog} The image of 
$\Phi_s:D^b_{coh}(C)\longrightarrow D^b_{coh}(X)$ is left orthogonal to the exceptional sequence $\sigma=
\Bigl({\cal O}_X(-d+3),\ldots,{\cal O}_X\Bigr) $ on $X$.
\eth
\pr
First, the sequence $\Bigl({\cal O}_X(-d+3),\ldots,{\cal O}_X\Bigr)$ is exceptional on $X$. Indeed, from the short exact sequence on a non-degenerate quadric $q$:
$$0\longrightarrow {\cal O}_q(i-2)
\longrightarrow {\cal O}_q(i)\longrightarrow {\cal O}_X(i)
\longrightarrow 0,$$
and from exceptionality of (\ref{qcol}) one can easily find: ${\rm H}^j\Bigl(X,{\cal O}(k)\Bigr)=0$, for any $j$ and $-d+3<k<0$.

Similarly to the proof of lemma 1.2 one can show the existence of the right adjoint functor $\Psi:D^b_{coh}(X)\longrightarrow D^b_{coh}(C)$ to $\Phi_S$.
Then for any object $A\in D^b_{coh}(C)$ one has
$${\rm Hom}_X\Bigl(\Phi_SA,{\cal O}(i)\Bigr)={\rm Hom}\Bigl(A,\Psi({\cal O}(i))\Bigr).$$
We have to show that $\Psi({\cal O}(i))=0$, for ${\cal O}(i)\in \sigma.$
 
Since there are no non-zero objects in $D^b_{coh}(C)$ which are orthogonal to all skyscraper sheaves ${\cal O}_c,\ c\in C$, it is sufficient to prove that all $\Phi_S({\cal O}_c)$ are orthogonal to $\sigma$. But every $\Phi_S({\cal O}_c)$ is a spinor bundle $S_q$ restricted to $X$ either from a smooth or from a degenerate quadric. In the former case we have an exact sequence
$$
0\longrightarrow S_q^{\vee}(-2)\longrightarrow S_q^{\vee}
\longrightarrow  S_q^{\vee}\left|_X\right.\longrightarrow 0.$$
Using this sequence and exceptionality of (\ref{qcol}) we easily find that $S_q$ is right orthogonal to $\sigma$.

If the quadric is degenerate, then we have a projection
$\pi:X\longrightarrow q'$ to a quadric of dimension $d-1$, and $\Phi({\cal O}_c)=\pi^*(S_{q'}).$ Therefore:
$${\rm Ext}^j_X\Bigl(\Phi({\cal O}_c),{\cal O}(i) \Bigr)=
{\rm Ext}^j_X\Bigl(\pi^*(S_{q'}),{\cal O}(i) \Bigr)=$$ $$={\rm Ext}^j_{q'}\Bigl(S_{q'},\pi_*{\cal O}(i) \Bigr)={\rm Ext}^j_{q'}\Bigl(S_{q'},{\cal O}(i)\oplus {\cal O}(i-1) \Bigr).$$
Because of exceptionality of (\ref{qcol2}) we are done.

Now we regard $\db{C}$ as a subcategory in  $\db{X}$.
As has been shown $\sigma=
\Bigl({\cal O}_X(-d+3),\ldots,{\cal O}_X\Bigr) $ lies in the right orthogonal to $\db{X}$.

\th{Theorem}\label{qdq}
The category $\db{X}$ on the intersection of two quadrics of dimension $d$ is generated as a triangulated category by
$\sigma $ and $\db{C}$, in other words there is a semiorthogonal decomposition
$$
\db{X}\ =\ \Bigl\langle{\cal O}_X(-d+3),\ldots,{\cal O}_X,\db{C}\Bigl\rangle
$$
\eth
\pr
Consider the subcategory $D\subset \db{X}$, generated by 
$\sigma $ and $\db{C}$. First, let us mention that the composition of the natural embedding $K_0(D)\otimes k\to
  K_0(X)\otimes k$ with the Chern character
$$
ch:K_0(X)\otimes k\to H^{even}(X,k)=\oplus H^{i,i}(X,k)
$$
is a surjective homomorphism from the Grothendieck group
of $D$, tensored by $k$, to the sum of the diagonal cohomologies of $X$ with coefficients in $k$. Indeed, the Chern character is a surjective morphism with the kernel, consisting, by Riemann-Roch-Hirzebruch formula, of those $v\in K_0(X)\otimes k$, which are in the (say, right) kernel
 of the Euler characteristic bilinear form $\chi $ (see ch.1).
From the orthogonal decomposition for $D$ one easily finds that the restriction of $\chi $ on $D$ has rank $d$, which coincides with the dimension of $H^*(X,k)$. The surjectivity
 follows.

Since $D$ has a semiorthogonal decomposition by admissible subcategories, it is in turn admissible \cite{BK}. Then, as usually, we suppose that $D^{\perp }$ is not trivial and consider an object $Z\in D^{\perp }$. It follows from above that $ch(Z)=0$.

Consider a singular quadric containing $X$. Let $\pi :X\to q$ be a projection from the singular point of this quadric to a non-singular quadric $q$ of dimension $d-1$.
There is a semiorthogonal decomposition
\begin{equation}\label{dqp}
\db{q}\ =\ \Bigl\langle{\cal O}_q(-d+2),\ldots,{\cal O}_q,S_q\Bigl\rangle
\end{equation}

For any $A\subset \db{q}$ we have an isomorphism:
$$
{\rm Hom}_X(\pi ^{*}A,Z)={\rm Hom}_q(A,\pi _{*}Z)
$$

Since all but the first element of (\ref{dqp}) after lifting to $X$ belong to the subcategory $D$, it follows that $\pi _{*}Z$
belongs to $\Bigl\langle{\cal O}_q(-d+3),\ldots,{\cal O}_q,S_q\Bigl\rangle ^{\perp }=\Bigl\langle{\cal O}_q(-d+2)\Bigl\rangle$.

We aim to prove that $Z$ as an object in $\db{X}$ is quasi-isomorphic to the direct sum of its cohomology sheaves.

\th{Lemma}
For any couple of coherent sheaves $A,B$ on $X$, such that 
$\pi_{*}A$ and $\pi_{*}B$ are direct sums of copies of a single linear bundle on $q$, one has
$$
{\rm Ext}^{i}_{X}(A,B)=0,\  {\rm for}\  i>1.
$$
\eth
\pr 
Let $s$ be the $\pi $-fibrewise involution on $X$. The fibred square $X'=X\times _qX$ of $X$ over $q$ is a union of two
copies of $X$, which normally intersects in the (smooth) ramification divisor $H$ of $\pi$ in $X$. These are the diagonal 
$$
\Delta X=\{(x,x)|x\in X\}
$$
and the $s$-diagonal
$$
\Delta _sX=\{(x,sx)|x\in X\}.
$$
This description implies a short exact sequence of coherent 
sheaves on $X'$:
\begin{equation}\label{diadec}
0\to {\cal O}_{\Delta _sX}(-H)\to {\cal O}_{X'}\to {\cal O}_{\Delta _X}\to 0.
\end{equation}

Denote by $p_{1}, p_{2}$ the projections of $X'$ to $X$. 
Take any coherent sheaf $C$ on $X$. Tensoring (\ref{diadec})
with $p_{1}^{*}C$ one obtains:
$$
0\to p_{1}^{*}C\otimes {\cal O}_{\Delta _sX}(-H)\to p_{1}^{*}C\to p_{1}^{*}C\otimes {\cal O}_{\Delta _X}\to 0
$$
Then, applying $p_{2*}$ to this sequence, one has:
$$
0\to s_{*}C(-1)\to p_{2*}p_{1}^{*}C\to C\to 0.
$$
Using the flat base change theorem (see \cite{H},II,\S5,prop.5.12) where the morphism and the base change both are $\pi $, we obtain:
$$
\pi ^{*}\pi _{*}=p_{2*}p_{1}^{*}
$$
Therefore one has an exact sequence for any sheaf on $X$:
\begin{equation}\label{ppdec}
0\to s_{*}C(-1)\to \pi ^{*}\pi _{*}C\to C\to 0.
\end{equation}

Let now $A$ and $B$ be such that $\pi_{*}A$ and $\pi_{*}B$
are sums of copies of a linear bundle on $q$. Without loss of generality we can assume that this linear bundle is trivial.
Since $\pi _{*}s_{*}\simeq \pi _{*}$, putting $C=B$ and $C=s_{*}B(-1)$ in (\ref{ppdec}), we obtain exact sequences:
$$
0\to s_{*}B(-1)\to \oplus {\cal O}\to B\to 0
$$
$$
0\to B(-2)\to \oplus {\cal O}(-1)\to s_{*}B(-1)\to 0.
$$

Juxtaposing these two sequences and then repeating the procedure in the same way one obtains a resolution for $B$:
$$
\dots \to \oplus {\cal O}(-2)\to \oplus {\cal O}(-1)\to \oplus {\cal O}\to B\to 0.
$$

Using this resolution and the fact that $\pi _{*}A$ is a trivial bundle, one obtains a spectral sequence converging to ${\rm Ext}^{\cdot }(B,A(-d+2))$:
$$
E^{p,q}_1= {\rm Ext}^q({\cal O}(-p)\:,\; {\cal O}(-d+2)) \Longrightarrow {\rm Ext}^{p+q}( B\:, \; A(-d+2))
$$
which shows that 
$$
{\rm Ext}^{<d-2}(B,A(-d+2))=0.
$$

By the adjunction formula one easily calculates the canonical class of $X$:
$$
\omega _X= {\cal O}_X(-d+2).
$$
Thus, by Serre duality the equality holds:
$$
{\rm Ext}^{i}( A\:, \; B)={\rm Ext}^{d-1-i}( B\:, \; A(-d+2))^{*}.
$$
It follows that ${\rm Ext}^{>1}( A\:, \; B)=0$.

Since ${\cal O}_q(-d+2)$ is an exceptional sheaf, any object from $\Bigl\langle{\cal O}_q(-d+2)\Bigl\rangle$ is isomorphic to the direct sum of its cohomology sheaves, which in turn are direct sums of copies of ${\cal O}_q(-d+2)$.

Let ${\cal H}^i$ be the cohomology sheaves of $Z$. Since $R^{0}\pi_{*}$ is an exact functor, $R^{0}\pi_{*}{\cal H}^i$
are direct sums of copies of ${\cal O}_q(-d+2)$. It follows by the lemma that
$$
{\rm Ext}^k({\cal H}^i\: ,\;{\cal H}^j)=0,\ {\rm for} k>1.
$$
It is well known that this implies a decomposition of $Z$ into a direct sum of its cohomology sheaves.

Since $Z\in D^{\perp }$, it follows that ${\cal H}^i\in D^{\perp }$. Then we have from above that $ch({\cal H}^i)=0$.
But a sheaf (not a complex of sheaves) with trivial Chern character is zero. As all cohomologies of $Z$ are zero, then $Z$ is quasi-isomorphic to zero itself. This finishes the proof of the theorem.

{\section{Birational transformations.}}

The aim of this chapter is to trace behaviour of the derived category of coherent sheaves with respect to birational transformations. It turns out that  blowing up and flip
transformations have a categorical  incarnation as adding or removing of semiorthogonal summands of a quite simple nature.

For a flop there are no such summands, thus it produces an equivalence of triangulated categories.

The simplest instance of a birational transformation is a blowing up of a variety along a smooth centre. A discription of the derived category of the blow-up in  terms of the categories of the variety and of the centre is done in \cite{or} . We give here the treatment with small corrections and with a stress on theorem \ref{mai} in proofs.

Let $Y$ be a smooth subvariety of codimension $r$ in a smooth algebraic variety $X$. Denote $\widetilde X$ the smooth algebraic variety obtained by the blowing up of $X$ along the centre $Y$.
There exists a fibred square:
$$
\begin{array}{ccc}
\widetilde Y&\stackrel{j }{\longrightarrow}& \widetilde X\\
\llap{\footnotesize $p$} \downarrow &&\llap{$\pi$} \downarrow\\
Y&\stackrel{i}{\longrightarrow}& X 
\end{array}
$$
where $i$ and $j$ are embeddings of smooth varieties, and $p: \widetilde Y\to Y$ is the projective fibration of the exceptional divisor $\widetilde Y$ in $\widetilde X$ over the centre $Y$. Recall that $\widetilde Y={\bf P}(N_{X/Y})$ is the projectivization of the normal bundle to $Y$ in $X$.

\th{Proposition}\label{blow} {\rm (see \cite{or})}
The pull back functors
$$
\pi^*: \db {X}\lto \db {\widetilde X}
$$
$$
p^*: \db {Y}\lto \db {\widetilde Y}
$$
are full and faithful.
\eth
\pr
The functor $\pi^*$ (resp., $p^*$) is isomorphic to $\Phi_{E_X}$ (resp., $\Phi_{E_Y}$), where $E_X$ (resp., $E_Y$) is the structure sheaf of the incidence subscheme $Z_X$ (resp., $Z_Y$) in $X\times {\widetilde X}$ (resp., $Y\times {\widetilde Y}$). 

We shall show that proof easily follows from theorem \ref{mai}. Indeed, for a point $x\in X$, $\Phi_{E_X}({\o x})=\pi^*({\o x})$ (resp., for a point $y\in Y$, $\Phi_{E_Y}({\o y})=p^*({\o y})$) is the structure sheaf of the $\pi$--fibre over $x$. Since fibres are disjoint,  orthogonality condition i) of theorem \ref{mai} follows. For $\Phi_{E_Y}$ analogously.

Since $\widetilde X$ has the same dimension as $X$ and due to the fact that for any couple $(F, E)$ of sheaves on a smooth variety ${\rm Ext}^i(F, E)=0$ for $i$ greater than the dimension of the variety, condition ii) of theorem \ref{mai} for $\Phi_{E_X}$ is verified.

Let us remark that for $A\subset B$, a smooth subvariety in a smooth variety, local extension groups for the structure sheaf ${\o A}$ of the subvariety are equal:
\begin{equation}\label{ext}
{\cal E}xt^i_B({\o A}\:, {\o A})=\Lambda^i N_{B/A},
\end{equation}
where $N_{B/A}$ is the normal vector bundle of $A$ in $B$.

For fibres $p^*({\o Y})$, which are biregulary isomorphic to projective spaces, the normal bundle is trivial. Thus they have no the higher cohomologies and the spectral sequence from local to global extension groups gives condition ii) of theorem \ref{mai} for $\Phi_{E_Y}$.

Another proof uses the projection formula:
$$
{\rm Ext}^i(p^*F\:,\:p^*G)={\rm Ext}^i(F\:,\:p_*p^*G)={\rm Ext}^i(F\:,\:p_*{\o {\widetilde Y}}\otimes G).
$$
Analogously for $\pi^*$.
Combining with the facts that $\pi_*{\o{\widetilde X}}={\o X}$ and $p_*{\o{\widetilde Y}}={\o Y}$ this gives the proof.
\th{Proposition}\label{emb}
For any invertible sheaf $L$ over $\widetilde Y$, the functor
$$
\Phi j_*(L\otimes p^*(\cdot)): \db {Y}\lto \db {\widetilde X}
$$
is full and faithful. 
\eth
\pr

Let us verify the hypothesis i) and ii) of the theorem \ref{mai}. For a point $y\in Y$ the image $\Phi ({\cal O}_y)$  is the structure sheaf of the corresponding $p$-fibre (regarded as a sheaf on ${\widetilde X}$). Since $p$-fibres over distinct points are disjoint the orthogonality condition    i) of the theorem \ref{mai} follows.

Now let us consider the structure sheaf ${\cal O}_F$ of a single $p$-fibre $F\subset {\widetilde X}$. By formula (\ref{ext}) one has the spectral sequence:
$$
{\rm E}^{i,j}_2=\oplus \ {\rm H}^i(F, \Lambda ^jN_{{\widetilde X}/F})=\oplus \ {\rm H}^i({\widetilde X},{\cal E}xt^j_{\widetilde X}({\o F}\:, {\o F}))\Longrightarrow \oplus \ 
{\rm Ext}^{i+j}_{\widetilde X}({\o F}\:, {\o F}).
$$

For the normal bundle $N_{{\widetilde X}/F}$ one has the exact sequence:
$$
0\to N_{{\widetilde Y}/F}\to N_{{\widetilde X}/F}\to N_{{\widetilde X}/{\widetilde Y}}|_F\to 0.
$$

Obviously, the normal bundle $N_{{\widetilde Y}/F}$ is trivial, while
$$
N_{{\widetilde X}/{\widetilde Y}}|_F\ =\ {\cal O}(Y)|_F\ =\ {\cal O}(-1).
$$
Since $F$ is a projective space and since there are no extension groups between ${\cal O}(-1)$ and ${\cal O}$ on a projective space, it follows that the above short exact sequence splits. Then the spectral sequence gives conditions ii) of the theorem \ref{mai}. This proves the proposition.

\th{Lemma}\label{combar}
Let $j:D\to Z$ be the embedding of a smooth divisor in a smooth algebraic variety $Z$ of dimension $n$. Consider for an object $A\in \db{X}$ an exact triangle with the canonical second morphism:
$$
\overline A\to j^{*}j_{*}A\to A
$$
Then $\overline A\simeq A\otimes {\cal O}_D(-D)[1]$.
\eth
\pr
The functor $j_{*}$ coincides with $\Phi _E$ (in notations of ch.1), where $E={\o G}$ is the structure sheaf of the graph subvariety $G$ for $j$ in $D\times Z$. The adjoint functor $j^{*}$ is isomorphic to $\Psi _{E'}$, where  $E'={\o G}$.

By proposition \ref{comp} in order to calculate $j^{*}j_{*}$ one has to find ${p_{13}}_*({p_{23}}^*{\o G}\otimes {p_{12}}^*{\o G})$, where $p_{ij}$ are the projections from the product $D\times Z\times D$ along the $k$-th component, where $\{ ijk\} =\{ 123\}$.

Note that
$$
{p_{12}}^*{\o G}={\o {G\times D}};\ \ {p_{23}}^*{\o G}={\o {D\times G}},
$$
where ${G\times D}$ and ${D\times G}$ are regarded as subvarieties in $D\times Z\times D$ of codimension $n$. These varieties intersect along the subvariety, which is the image of the morphism $({\rm id},j,{\rm id}):D\to D\times Z\times D$. It is of codimension $2n-1$, hence a non-transversal intersection.

Fortunately, both ${G\times D}$ and ${D\times G}$ lie in the image of $\Delta_3:D^3\to D\times Z\times D$, where they meet transversally.

This helps to compute ${\cal H}^i({\o {G\times D}}\otimes {\o {D\times G}})$. Indeed, one can consider the tensor product ${\o {G\times D}}\otimes {\o {D\times G}}$ as the restriction of ${\o {D\times G}}$ to ${G\times D}$. Restricting first to the divisor $D^3$, we obtain:
$$
\begin{array}{l}
L^0\Delta_3{\o {D\times G}}={\o {D\times G}},\\ L^1\Delta_3{\o {D\times G}}={\o {D\times G}}(-D^3),\\
L^i\Delta_3{\o {D\times G}}=0, \ {\rm for}\  i>1.
\end{array}
$$

Then restriction to  ${G\times D}\subset D\times Z\times D$ and projection along $Z$ give that the complex ${\cal K}={p_{13}}_*({p_{23}}^*{\o G}\otimes {p_{12}}^*{\o G})\in {\db {D\times Z}}$ has only two cohomology sheaves. Namely,
$$
\begin{array}{l}
{\cal H}^0({\cal K})={\o {\Delta}},\\
{\cal H}^{-1}({\cal K})={\o {\Delta }}\otimes \pi ^*{\cal O}(-D).
\end{array}
$$

Therefore one has the exact triangle
$$
{\o {\Delta }}\otimes \pi ^*{\cal O}(-D)[1]\to {\cal K}\to {\o {\Delta }}.
$$
Applying functors, corresponding to objects from this triangle to arbitrary $A\in {\db  D}$ one gets the proof.

Denote by $D(X)$ the full subcategory of $\db{\widetilde X}$ which is the image of $\db{X}$ with respect to the functor $\pi^*$ and by $D(Y)_k$ the full subcategories of $\db{\widetilde X}$ which are the images of $\db{Y}$ with respect to the functors $j_*(\o{\widetilde Y}(k)\otimes p^*(\cdot))$.
 

\th{Proposition}\label{ort}
The sequence
$$
\Bigl( D(Y)_{-r+1}, ... , D(Y)_{-1}, D(X) \Bigr)
$$
is a semiorthogonal sequence of admissible subcategories in $\db{\widetilde X}$.
\eth
\pr
1). Let $j_*A\in D(Y)_k$ and $j_*B\in D(Y)_m$ with $r-2\ge k-m >0$ .

That means
\begin{equation}\label{fb}
A=p^*A'\otimes \o{\widetilde Y}(k),\; B=p^*B'\otimes \o{\widetilde Y}(m),
\end{equation}
for some  $A', B'\in \db{Y}$.

We have an exact triangle:
\begin{equation}\label{jj}
\bar A \lto j^*j_*A\lto A,
\end{equation}
and by lemma \ref{combar} an isomorphism 
$$
\bar A\cong A\otimes \o{\widetilde Y}(1)[1].
$$

Furthermore,
$$
\begin{array}{l}
{\h A, B}\cong {\h p^*A'\otimes \o{\widetilde Y}(k), {p^*B'\otimes \o{\widetilde Y}(m)}}\cong
\\
{\h p^*A', {p^*B'\otimes \o{\widetilde Y}(m-k)}}\cong {\h A', {p_*(p^*B'\otimes \o{\widetilde Y}(m-k))}}\cong\\
{\h A', {B'\otimes p_*\o{\widetilde Y}(m-k)}}.
\end{array}
$$
From  vanishing of $p_*\o{\widetilde Y}(-n)=0$, with $r-1\ge n >0$ we obtain:
$$
{\h A, B}=0.
$$
Analogously
$$
{\h {\bar A}, B}=0.
$$
Then, triangle (\ref{jj}) gives
$$
{\h j_*A, {j_*B}}=0.
$$
This proves semiorthogonality for the sequence of the subcategories
$$
\Bigl( D(Y)_{-r+1}, ... , D(Y)_{-1} \Bigl).
$$
2). If $\pi^*A\in D(X)$ and $j_*B\in D(Y)_m$ for $-r+1\le m\le -1$ with $B$ being of the form (\ref{fb}), then
$$
{\h \pi^*A, {j_*B}}\cong {\h A, {\pi_*j_*B}}\cong {\h A, {i_*p_*B}}=0.
$$
This is equal to zero, because $p_*B=p_*(p^*B'\otimes \o{\widetilde Y}(m))=B'\otimes p_*\o{\widetilde Y}(m)=0$.

\th{Theorem}\label{full}
In the above notations, the semiorthogonal sequence of admissible subcategories
$$
\Bigl< D(Y)_{-r+1}, ... , D(Y)_{-1}, D(X) \Bigl>
$$
generates  the category $\db{\widetilde X}$.
\eth
\pr
See \cite{or}.

This theorem gives a semiorthogonal decomposition of the derived category $\db{\widetilde X}$ on a blow-up $\widetilde X$. It was used in \cite{or} for constructing a full exceptional sequence in $\db{\widetilde X}$, starting from such sequences on $X$ and $Y$.
\bigskip

Now we explore the behaviour of the derived categories of coherent sheaves with respect to simplest flip and flop transformations.

Let $Y$ be a smooth subvariety of a smooth algebraic variety
$X$ such that $Y\cong {\bf P}^k$ and $N_{X/Y}\cong {\o Y}(-1)^{\oplus (l+1)}$ with $l\le k$.

If now $\widetilde X$ is a blow-up of $X$ along $Y$, then  exceptional divisor $\widetilde Y\cong {\bf P}^k\times {\bf P}^l$ is isomorphic to the product of projective spaces. This allows us to blow down $\widetilde X$ in such a way that$\widetilde Y$ project to the second component ${\bf P}^l$ of the product. As a result we obtain a smooth algebraic variety $X^+$ with subvariety $Y^+\cong {\bf P}^l$. This situation is depicted in the following diagram:

\begin{figure}[h]
\hspace*{6cm}\epsffile{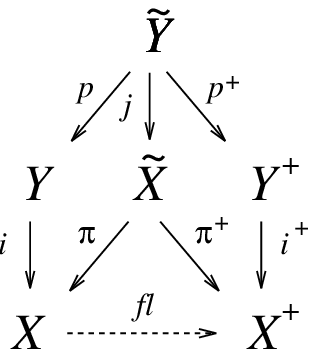}
\end{figure}

 The birational map $X\lto X^+$ is a simple example of flip, for $l<k$, and flop, for $l=k$, transformations.

One can easily calculate that for the restriction $\o{\widetilde X}(\widetilde Y){\Bigl |}_{\widetilde Y}$ there exists an isomorphism
$$
{\o {\widetilde X}}(\widetilde Y){\Bigl |}_{\widetilde Y}\cong {\cal O}(-1)\b {\cal O}(-1),
$$
where ${\cal O}(-1)\b {\cal O}(-1):= p^*{\o Y}(-1)\otimes p^{+*}{\o {Y^+}}(-1)$.
For subsequent calculations we need the formula for the canonical sheaf $\omega_{\widetilde X}$ of the blow-up.
$$
\omega_{\widetilde X}\cong \pi^*\omega_X\otimes \o{\widetilde X}(l\widetilde Y).
$$
Further, by the adjunction formula we know that
$$
\omega_X{\Bigl |}_Y\cong \omega_Y\otimes \Lambda^{l+1}N^*_{X/Y}\cong {\o Y}(l-k).
$$
Combining these facts we conclude that
$$
\omega_{\widetilde X}{\Bigl |}_{\widetilde Y}\cong (\pi^*\omega_X\otimes {\o {\widetilde X}}(l\widetilde Y)){\Bigl |}_{\widetilde Y}\cong p^*(\omega_X\Bigl|_Y)\otimes {\o {\widetilde X}}(l\widetilde Y){\Bigl |}_{\widetilde Y}\cong {\cal O}(-k)\b {\cal O}(-l).
$$
Next is the main theorem of this section.


\th{Theorem}\label{fl}
In the above notations, the functor
$$
\pi_*\pi^{+*}: \db {X^+}\lto \db {X}
$$
is full and faithful.
\eth
\pr
 We have to show that for any pair of objects $A, B\subset \db {X^+}$ there is an isomorphism
$$
{\h \pi_*\pi^{+*}A, {\pi_*\pi^{+*}B}}\cong {\h A, B}.
$$
For the left hand side we have a canonical isomorphism 
\begin{equation}\label{pi}
{\h \pi_*\pi^{+*}A, {\pi_*\pi^{+*}B}}\cong {\h \pi^* \pi_*\pi^{+*}A, {\pi^{+*}B}}.
\end{equation}
Consider an exact triangle:
\begin{equation}\label{nbar}
 \pi^* \pi_*\pi^{+*}A\lto \pi^{+*}A\lto {\bar A}. 
\end{equation}
Applying to it the functor ${\h {\cdot}, {\pi^{+*}B}}$ we find  that if
\begin{equation}\label{zer}
{\h {\bar A}, {\pi^{+*}B}}=0, 
\end{equation}
then the right hand side of (\ref{pi}) is isomorphic to ${\h \pi^{+*}A, {\pi^{+*}B}}$. Since $\pi^+$ is an instance of a blow up morphism, by proposition \ref{blow} we have an isomorphism:
$$
{\h \pi^{+*}A, {\pi^{+*}B}}\cong {\h A, B}.
$$
Hence we need to check (\ref{zer}).

Again by proposition \ref{blow} $\pi^*$ is full and faithful, or, in other words, $\pi_*\pi^*$ is isomorphic to the identity functor. Applying this isomorphism to the object $\pi_*\pi^{+*}A$, we get the first morphism from the exact triangle, obtained by application of functor $\pi_*$ to triangle (\ref{nbar}):
$$
\begin{array}{ccccc}
\pi_*\pi^*\pi_*\pi^{+*}A&\stackrel{\sim}{ \lto}&\pi_*\pi^{+*}A&\lto&\pi_*{\bar A}.
\end{array}
$$

Therefore, $\pi_*{\bar A}=0$. Consequently, for any object $C\in\db {X^+}$ ${\h \pi^*C , {\bar A}}=0$, i.e. $\bar A\in D(X)^{\perp}$.

Recall that by  theorem \ref{full}
$$
D(X)^{\perp}=\Bigl< D(Y)_{-l}, ... , D(Y)_{-1} \Bigl>
$$
is a semiorthogonal decomposition of the category $D(X)^{\perp}$. The notations $D(Y)_{-k}$ are fixed before proposition \ref{ort}.

If we choose full exceptional sequences in each $D(Y)_{-i}$, then gathering them together we obtain a full sequence in $D(X)^{\perp}$. The following one will be convinient for us:
$$
\begin{array}{lllll}
D(X)^{\perp}=\Bigl< &j_*({\cal O}(-k)\b {\cal O}(-l)),& ... & ....& j_*({\cal O}\b {\cal O}(-l)),\\
& j_*({\cal O}(-k+1)\b {\cal O}(-l+1)),& ... & ....
& j_*({\cal O}(1)\b {\cal O}(-l+1)),\\
& ................ & ... & ... & .............. \\
& j_*({\cal O}(-k+l-1)\b {\cal O}(-1)),&... & ...& j_*({\cal O}(-l-1)\b {\cal O}(-1))\Bigl>
\end{array}
$$

Let us divide this sequence in two parts ${\cal A}$ and ${\cal B}$, such that
$$
D(X)^{\perp}=\Bigl< {\cal B , A}\Bigl>
$$
be a semiorthogonal decomposition for $D(X)^{\perp}$
with ${\cal A}$ and ${\cal B}$ being the subcategories generated by $j_*({\cal O}(i)\b {\cal O}(s))$ with $i\ge 0$ and $i<0$ respectively.
If $1\le i\le k$ and $1\le s\le l$ then the object
$ j_*({\cal O}(-i)\b {\cal O}(-s))$ belongs to simultaneously $D(X)^{\perp}$ and $D(X^+)^{\perp}$. Therefore, applying the functor Hom with the target in this object to exact triangle (\ref{nbar})
we obtain:
$$
{\h {\bar A}, { j_*({\cal O}(-i)\b {\cal O}(-s))}}=0,\qquad \mbox{for } 1\le i\le k \;\mbox{and } 1\le j \le l.
$$

Since ${\bar A}\in D(X)^{\perp}$, it immediately follows that $\bar A\in {\cal A}$.
Further, we notice that ${\cal A}\otimes\omega_{\widetilde X}\in D(X^+)^{\perp}$, because $\omega_{\widetilde X}{\Bigl |}_{\widetilde Y}\cong {\cal O}(-k)\b {\cal O}(-l)$ and $l\le k$. Therefore, for $B\in{\db {X^+}}$
$$
{\h \pi^{+*}B , {\bar A \otimes \omega_{\widetilde X}}}=0.
$$
Hence by Serre duality (\ref{zer}) follows.
This proves the theorem.
\bigskip

{\sc Remark}. For the case of flop ($l=k$) the functor $\pi_*\pi^{+*}$ is an equivalence of triangulated categories.
\bigskip


Now we investigate more carefully 3-dimensional flops.

Let $f:X\longrightarrow Y$ be a proper birational morphism between compact threefolds, which blow down only an indecomposable curve $C$. Assume that $X$ is smooth and $C\cdot K_X=0$. Then $C\cong{\bf CP }^1$ and $N_{X/C}$ is equal to either ${\cal O}(-1)\oplus{\cal O}(-1)$ or 
 ${\cal O}\oplus{\cal O}(-2)$, or  ${\cal O}(1)\oplus{\cal O}(-3)$ (see, e.g. \cite{CKM}). There exist in this situation a smooth compact threefold $X^+$ with a curve $C^+\subset X^+$ and with a morphism $f^+:X^+\longrightarrow Y$, which blows down only the curve $C^+$, and with birational, but not biregular, map $g:X\longrightarrow X^+$, which is embedded in the commutative triangle
\begin{figure}[h]
\hspace*{6cm}\epsffile{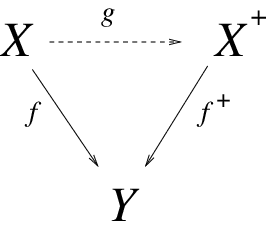}
\end{figure}


Such $X^+$ is unique (see \cite{K}). Birational map $g$ is called flop; $g$ induces isomorphism $X\setminus C\stackrel{\sim}{\rightarrow}X^+\setminus C^+$. Let us remark that the curve $C^+$ also has trivial intersection with $K_{X^+}$ and its normal bundle is of the same kind as the one on $C$.

If the curve $C$ has the normal bundle either of the first or of the second kind, i.e., isomorphic to ${\cal O}(-1)\oplus{\cal O}(-1)$ or 
 ${\cal O}\oplus{\cal O}(-2)$, then, following M.~Reid \cite{R}, we call it $(-2)$-curve.

We are going to prove that if $X$ and $X^+$ are related by flop with $(-2)$-curve, then $D^b_{coh}(X)$ is equivalent to
$D^b_{coh}(X^+)$. 

This supplies the following

{\bf Conjecture.} If two smooth varieties are related by flop, then the derived categories of coherent sheaves on them are equivalent as triangulated categories.

{\bf Comment.} By \cite{K} any birational transformation between two 3-dimensional Calabi--Yau varieties can be decomposed in a sequence of flops. Therefore the conjecture would imply equivalence of the derived categories of any two birationally isomorphic 3-dimensional Calabi--Yau's.

\vspace{2ex}

Let us remark that there exist examples of flops on threefolds with the normal bundle $N_{X/C}={\cal O}(1)\oplus {\cal O}(-3)$.

Consider a smooth compact threefold $X$ with a curve $C\cong {\bf CP}^1$, which is a $(-2)$-curve. Then there exist a flop
$X\longrightarrow X^+$. It is known an explicit decomposition for it in the so called `pagoda' of M.~Reid:

\begin{figure}[h]
\hspace*{6cm}\epsffile{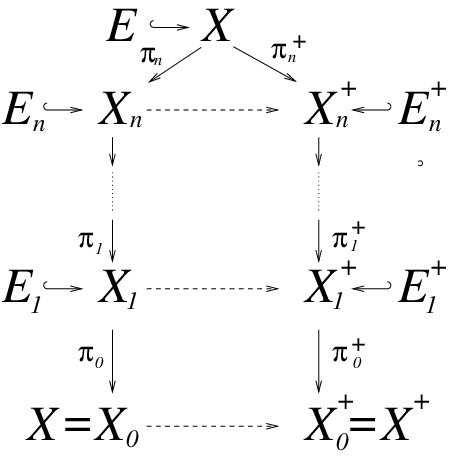}
\end{figure}

Here $X_1$ is a blow-up of $X$ in a curve $C$, $E_1$ the  exceptional divisor of this transformation, which is isomorphic to ${\bf F}_2$. The exceptional section $S_1\hookrightarrow E_1$ also is a $(-2)$-curve on $X_1$, so $X_2$ is the blow-up of $X_1$ in $S_1$ and so on. Finally we obtain a threefold $X_n$ with a divisor $E_n\simeq {\bf F}_2$ and with the section $S_n$, such that $N_{X_n/S_n}={\cal O}(-1)\oplus {\cal O}(-1)$. The blow-up of $X_n$ in $S_n$ is a threefold $\widetilde X$ with the exceptional divisor $E\simeq 
{\bf P}^1\times {\bf P}^1$. Contracting it in the other direction we obtain $X^+_n$. Further, contracting one-by-one
the proper transforms of divisors $E_n,\ldots,E_1$ we break our way through to a threefold $X^+$. See detals in the original paper \cite{R}.

All maps $X_i\longrightarrow X_i^+$ are flops too. We are going to proceed by induction on the length of the `pagoda'.

Denote by $\pi_i$ the birational morphisms $X_{i+1}\longrightarrow X_i$ and by $\Pi_i$ the morphisms $\widetilde X\longrightarrow X_i$ obtained by composition. Similarly (by $\pi^+_i$ and $\Pi^+_i$ ) for the right side of the `pagoda'. By proposition \ref{blow} the pull back functor for a blowing up is full and faithful, therefore the derived categories $D^b_{coh}(X_i)$ and $D^b_{coh}(X^+_i)$ can be identified with the full admissible subcategories in $D^b_{coh}(\widetilde X)$. Denote them $D(X_i)$ and $D(X^+_i)$ respectively. We have two filtrations on $D^b_{coh}(\widetilde X)$:
$$D(X_0)\subset D(X_1)\subset\ldots \subset D(X_n)\subset
D^b_{coh}(\widetilde X),$$
$$D(X^+_0)\subset D(X^+_1)\subset\ldots \subset D(X^+_n)\subset
D^b_{coh}(\widetilde X).$$
Denote ${\cal A}_i$ (resp.,
${\cal A}^+_i$  ) the right orthogonal to $D(X_i)$ (resp.,
$D(X^+_i)$), i.e. ${\cal A}_i=D(X_i)^{\bot},\
{\cal A}^+_i=D(X^+_i)^{\bot}.$
Denote ${\cal B}_i$ the common part of ${\cal A}_i$ and
${\cal A}^+_i $, i.e. the full subcategory, consisting of the objects, which are right orthogonal to $D(X_i)$ and to
$D(X^+_i)$.

\th{Proposition}\label{asub}  In the above notations one has

i) ${\cal B}_i$ are admissible subcategories

ii) there exists  a semiorthogonal decomposition of the categories
${\cal A}_i$ and ${\cal A}^+_i $ in pairs of admissible subcategories
$${\cal A}_i=\Bigl\langle{\cal B}_i,{\cal C}_i\Bigr\rangle\qquad{\cal A}^+_i=\Bigl\langle{\cal D}_i,{\cal B}_i\Bigr\rangle ,$$
such that ${\cal D}_i={\cal C}_i\otimes\omega_{\widetilde X}$
(i.e. subcategory ${\cal D}_i$ consists of those objects, which are twists by the canonical sheaf $\omega_{\widetilde X}$
of the objects from ${\cal C}_i$).
\eth
\pr We use induction by the length of the `pagoda'. The base of the induction is a flop in a curve with normal bundle ${\cal O}(-1)\oplus{\cal O}(-1)$.

In our notations we have the subcategories $D(X_n)$ and $D(X^+_n)$. Moreover, we can choose the following decompositions for ${\cal A}_n$ and  ${\cal A}^+_n$
(theorem \ref{full})
$${\cal A}_n=\Bigl\langle {\cal O}_E(-l'-l''),{\cal O}_E(-l'')\Bigr\rangle ,$$
$${\cal A}^+_n=\Bigl\langle {\cal O}_E(-l'-2l''),{\cal O}_E(-l'-l'')\Bigr\rangle ,$$
where $l'$ and $l''$ are fibres of the projections of $E$ to $S_n$ and $S^+_n$ respectively. It follows that ${\cal B}_n$ is a subcategory generated by one exceptional object ${\cal O}_E(-l'-l'')$, hence admissible, and ${\cal C}_n$ and
${\cal D}_n$ are also generated by one object ${\cal O}_E(-l'')$ and  ${\cal O}_E(-l'2-l'')$  respectively.

We have the formula for the restriction of $\omega_{\widetilde X}$ to $E$:
$$\omega_{\widetilde X}\left|_E\right.={\cal O}(-l'-l'').$$
It follows that ${\cal D}_n={\cal C}_n\otimes\omega_{\widetilde X}.$

Now suppose that for $i>0$ we have already proved that ${\cal B}_i$ are admissible, and that ${\cal A}_i$ and ${\cal A}_i^+$ have semiorthogonal decompositions
$$
{\cal A}_i=\Bigl\langle {\cal B}_i,\; {\cal C}_i\Bigl\rangle\quad\mbox{and }\; {\cal A}_i^+=\Bigl\langle {\cal D}_i,\; {\cal B}_i\Bigl\rangle
$$
$$
\mbox{with}\quad{\cal D}_i={\cal C}_i\otimes\omega_{\tilde X}.
$$
Let us prove it for $i=0$.

Again by theorem \ref{full} we have the decomposition
$$
D(X_1)=\Bigl\langle \Pi^*_1{\o {E_1}}(-s_1-2l_1),\; \Pi^*_1{\o {E_1}}(-s_1-l_1), \;D(X_0)\Bigl\rangle,
$$
where $s_1$ is the class of the exceptional section of $E_1\simeq{\bf F}_2$ and $l_1$ is a fibre. The decomposition for ${\cal A}_0$ follows:
$$
{\cal A}_0=\Bigl\langle {\cal A}_1,\; \Pi^*_1{\o {E_1}}(-s_1-2l_1),\; \Pi^*_1{\o {E_1}}(-s_1-l_1)\Bigl\rangle.
$$
Now we shall show that
$$
\Pi^*_1{\o {E_1}}(-s_1-2l_1)=\Pi^{+*}_1{\o {E_1^+}}(-s_1^+-2l_1^+).
$$
There exists an exact sequence on $X_1$:
$$
0\lto {\o {E_1}}(-s_1-2l_1)\lto{\o {E_1}}\lto{\o \Gamma}\lto 0,
$$
where $\Gamma$ is a curve from the linear system $|s_1+2l_1|$ on $E_1$. The main point here is that $\Gamma\cap  S_1=\emptyset$, i.e. $\Gamma$ does not intersect with the locus for the blowing up of $X_1$. Therefore
$$
\Pi^*_1{\o {\Gamma}}\cong\Pi^{+*}_1{\o {\Gamma}}
$$
(we identify the curve $\Gamma$ with its proper transform on $X^+_1$).

Moreover, $\Pi^*_1{\o {E_1}}={\o Z}$, where $Z=\bigcup^n_{i=1}E_i\bigcup E$, i.e. ${\o Z}=\Pi^{+*}_1{\o {E_1^+}}$. Again using the fact that the pull back functor for a blowing up is full and faithful, we immediately obtain
\begin{equation}\label{=1}
\Pi^*_1{\o {E_1}}(-s_1-2l_1)=\Pi^{+*}_1{\o {E_1^+}}(-s_1^+-2l_1^+).
\end{equation}
Consider again the decomposition for ${\cal A}_0$:
$$
{\cal A}_0=\Bigl\langle {\cal A}_1,\; \Pi^*_1{\o {E_1}}(-s_1-2l_1),\; \Pi^*_1{\o {E_1}}(-s_1-l_1)\Bigl\rangle=
$$
$$
=\Bigl\langle {\cal B}_1,\; {\cal C}_1, \Pi^*_1{\o {E_1}}(-s_1-2l_1),\; \Pi^*_1{\o {E_1}}(-s_1-l_1)\Bigl\rangle.
$$
For any object $C\subset{\cal C}_1$ one has:
$$
{\h C, {\Pi^*_1{\o {E_1}}(-s_1-2l_1)}}\cong{\h \Pi^*_1{\o {E_1}}(-s_1-2l_1), {C\otimes\omega_{\tilde X}[3]}}^*\cong
$$
$$
\cong{\h \Pi^{+*}_1{\o {E_1^+}}(-s_1^+-2l_1^+), {C\otimes\omega_{\tilde X}[3]}}^*=0.
$$
The last equation is due to the fact that $C\otimes\omega_{\tilde X}\in{\cal D}_1$, i.e. it belongs to $D(X^+_1)^{\perp}$.

Therefore, the subcategory ${\cal C}_1$ and the object $ \Pi^*_1{\o {E_1}}(-s_1-2l_1)$ are both sides mutually orthogonal. This means that we can exchange their positions in the semiorthogonal decomposition
$$
{\cal A}_0=\Bigl\langle {\cal B}_1,\;  \Pi^*_1{\o {E_1}}(-s_1-2l_1),\; {\cal C}_1, \;\Pi^*_1{\o {E_1}}(-s_1-l_1)\Bigl\rangle .
$$
It follows from (\ref{=1}) that the object $\Pi^*_1{\o {E_1}}(-s_1-2l_1)$ is orthogonal from the right to both $D(X_0)$ and $D(X_0^+)$.

Therefore, we have the semiorthogonal decomposition for ${\cal B}_0$:
$$
{\cal B}_0=\Bigl\langle {\cal B}_1,\; \Pi^*_1{\o {E_1}}(-s_1-2l_1) \Bigl\rangle,
$$
hence ${\cal B}_0$ is admissible.

For ${\cal C}_0$ we have the decomposition:
$$
{\cal C}_0=\Bigl\langle {\cal C}_1,\; \Pi^*_1{\o {E_1}}(-s_1-l_1)\Bigl\rangle.
$$
Now let us consider a subcategory ${\cal A}_0^+$. We can choose the following decomposition for it
$$
{\cal A}_0^+=\Bigl\langle {\cal A}_1^+,\; \Pi^{+*}_1{\o {E_1^+}}(-s_1^+-3l_1^+),\; \Pi^{+*}_1{\o {E_1^+}}(-s_1^+-2l_1^+)\Bigl\rangle=
$$
$$
=\Bigl\langle {\cal D}_1,\; {\cal B}_1, \Pi^{+*}_1{\o {E_1^+}}(-s_1^+-3l_1^+),\; \Pi^{+*}_1{\o {E_1^+}}(-s_1^+-2l_1^+)\Bigl\rangle.
$$

Now we shall show that
\begin{equation}\label{=2}
\Pi^{+*}_1{\o {E_1^+}}(-s_1^+-3l_1^+)\cong \Pi^*_1{\o {E_1}}(-s_1-l_1)\otimes\omega_{\tilde X}.
\end{equation}
It follows from here that:

first, the objects ${\cal B}_1$ and $\Pi^{+*}_1{\o {E_1^+}}(-s_1^+-3l_1^+)$ are both sides mutually orthogonal.
Therefore, one can exchange their positions in the decomposition for ${\cal A}_0^+$, because for any object $B\in{\cal B}_1$
$$
{\h B, {\Pi^{+*}_1{\o {E_1^+}}(-s_1^+-3l_1^+)}}\cong{\h {\Pi^*_1{\o {E_1}}(-s_1-l_1)[-3]}, B}^*=0.
$$

second, for ${\cal D}_0$ one has the decomposition
$$ 
{\cal D}_0=\Bigl\langle {\cal D}_1,\;  \Pi^{+*}_1{\o {E_1^+}}(-s_1^+-3l_1^+)\Bigl\rangle.
$$
This implies that ${\cal D}_0={\cal C}_0\otimes\omega_{\tilde X}$, because for ${\cal D}_1$ we have such decomposition by induction, hence (\ref{=2}) allows us to claim this for ${\cal D}_0$.

Thus, the proof of the proposition follows from the
\th{Lemma}\label{congr}
In the above notations one has:
$$
\Pi^{+*}_1{\o {E_1^+}}(-s_1^+-3l_1^+)\cong \Pi^*_1{\o {E_1}}(-s_1-l_1)\otimes\omega_{\tilde X}.
$$
\eth
\pr
Standart calculations for blow-ups give that the restriction of $\omega_{ X_1}$ to $E_1$ is an invertible sheaf ${\o {E_1}}(-s_1-2l_1)$. Further, $\Pi_1^*\omega_{X_1}\cong\Pi_1^{+*}\omega_{X_1^+}$, because
$$
\omega_{\tilde X}\cong\Pi_1^*\omega_{X_1}\otimes{\o {\tilde X}}(2E_2+\cdots 2E_n +E)
$$
and
$$
 \omega_{\tilde X}\cong\Pi_1^{+*}\omega_{X_1^+}\otimes{\o {\tilde X}}(2E_2+\cdots 2E_n +E)
$$
(we denote by common letter $E_i$ the exceptional divisor on $X_i$ as well as its proper transforms on $X_{i+1},...,{\widetilde X}$).

One has
$$
{\o {E_1}}(-s_1-l_1)={\o {E_1}}(l_1)\otimes\omega_{X_1},
$$
$$
{\o {E_1^+}}(-s_1^+-3l_1^+)={\o {E_1^+}}(-l_1^+)\otimes\omega_{X_1^+}.
$$
Therefore, it is sufficient to show that
$$
\Pi^{+*}_1{\o {E_1^+}}(-l_1^+)\cong \Pi^*_1{\o {E_1}}(l_1)\otimes\omega_{\tilde X}.
$$
 
This fact we shall prove also by induction on the length of the `pagoda'.

The base of the induction:  Consider $X_n$ and ${\o {E_n}}(l_n)$. Then for $\Pi^*_n{\o {E_n}}(l_n)$ we have the exact sequence:
$$
0\lto{\o {E_n}}(-s_n+l_n)\lto\Pi^*_n{\o {E_n}}(l_n)\lto 
{\o E}(l')\lto 0.
$$
Let us twist it by $\omega_{\tilde X}$. We know that
$$
\omega_{\tilde X}\cong\Pi_n^*\omega_{X_n}\otimes{\o {\tilde X}}(E),
$$
and
$$
\omega_{X_n}\Bigl|_{E_n}={\o {E_n}}(-s_n-2l_n), {\o {\tilde X}}(E)\Bigl|_{E_n}={\o {E_n}}(s_n),
$$
$$
\omega_{\tilde X}={\o E}(-l'-l'').
$$
Therefore
$$
0\longrightarrow {\cal O}_{E_n}(-s_n-l_n)\longrightarrow \Pi^*_n {\cal O}_{E_n}(l_n)\otimes\omega_{\widetilde X}\longrightarrow{\cal O}_{E}(-l'')\longrightarrow 0.
$$
From the other hand, for $\Pi^{+*}_n {\cal O}_{E^+_n}(-l^+_n)$
we have the short exact sequence
$$
0\longrightarrow {\cal O}_{E^+_n}(-s^+_n-l^+_n)\longrightarrow \Pi^{+*}_n {\cal O}_{E^+_n}(-l^+_n)\longrightarrow{\cal O}_{E}(-l'')\longrightarrow 0.
$$
Keeping in mind that $E_n^+=E_n$ on $\widetilde X$ and that ${\rm Ext}^1({\cal O}_E(-l')\;,\; {\cal O}_E(-s_n-l_n))$ is one dimensional we conclude that
$$
\Pi_n^*{\cal O}_{E_n}(l_n)\otimes \omega_{\widetilde X}\cong\Pi_n^{+*}{\cal O}_{E^+_n}(-l^+_n).
$$
To make one step of the induction is in this case practically the same as to check the base of the induction. Namely, for 
$\pi_1^*{\cal O}_{E_1}(l_1)$ one has the short exact sequence on $X_2$:
$$
0\longrightarrow {\cal O}_{E_1}(-s_1+l_1)\longrightarrow \pi^{*}_1 {\cal O}_{E_1}(l_1)\longrightarrow{\cal O}_{E_2}(l_2)\longrightarrow 0.
$$
Lifting it up to ${\widetilde X}$ and twisting by $\omega_{\widetilde X}$ one obtains:
$$
0\longrightarrow {\cal O}_{E_1}(-s_1-l_1)\longrightarrow \Pi^*_1 {\cal O}_{E_1}(l_1)\otimes \omega_{\widetilde X}\longrightarrow\Pi^*_2{\cal O}_{E_2}(l_2)\otimes\omega_{\widetilde X}\longrightarrow 0.
$$
(we use here the fact that $E_1\cap E_2=\emptyset$, for $i>2$).

By hypothesis of the induction
$$
\Pi_2^*{\cal O}_{E_2}(l_2)\otimes\omega_{\widetilde X}\cong
\Pi_2^{+*}{\cal O}_{E^+_2}(-l_2),
$$
and the sheaf ${\cal O}_{E_1}(-s_1-l_1)$ coincides with ${\cal O}_{E_1^+}(-s_1^+-l_1^+)$.
Using as above that
${\rm Ext}^1(\Pi_2^*{\cal O}_{E_2^+}(-l_2)\;,\; {\cal O}_{E_1^+}(-s_1^+-l^+_1))$ is of dimension 1, we obtain:
$$
\Pi_1^*{\cal O}_{E_1}(l_1)\otimes\omega_{\widetilde X}\cong
\Pi_1^{+*}{\cal O}_{E^+_1}(-l_1).
$$
This proves the lemma and, consequently,  proposition \ref{asub}.

\th{Theorem}\label{efl}
The functor
$$
\Pi_*\Pi^{+*}:\db{X^+}\lto\db{X}
$$
is an equivalence of triangulated categories.
\eth
\pr
Let us consider two objects $A, B\in\db{X^+}$, then we have:
$$
{\h \Pi_*\Pi^{+*}A, {\Pi_*\Pi^{+*}B}}\cong{\h \Pi^*\Pi_*\Pi^{+*}A, {\Pi^{+*}B}}.
$$
There exists an exact triangle:
$$
\Pi^*\Pi_*\Pi^{+*}A\lto\Pi^{+*}A\lto{\bar A}.
$$

Since by proposition \ref{blow} $\Pi_*\Pi^*$ is isomorphic to the identity functor on $\db{X}$, one can easily see that ${\bar A}\in D(X)^{\perp}$.
Moreover the group of the homomorphisms from ${\bar A}$ to any object of the subcategory ${\cal B}_0$ is trivial, for so are the groups of the homomorphisms from the other members of the exact triangle.
Therefore ${\bar A}\in{\cal C}$.

It follows that 
$$
{\h {\bar A}, {\Pi^{+*}B}}\cong{\h {\Pi^{+*}B}, {{\bar A}\otimes\omega_{\tilde X}[3]}}^*=0,
$$
because ${\bar A}\otimes\omega_{\tilde X}\in{\cal D}\subset{D(X^+)}^{\perp}$.

Therefore,
$$
{\h \Pi_*\Pi^{+*}A, {\Pi_*\Pi^{+*}B}}\cong{\h \Pi^{+*}A, {\Pi^{+*}B}}\cong{\h A, B}.
$$
The latter isomorphism is due to fully faithfulness of $\Pi^{+*}$ (by proposition \ref{blow}).

This proves that $\Pi_*\Pi^{+*}$ is full and faithful.

Now suppose that $\Pi_*\Pi^{+*}$ is not an equivalence. Then $\Pi_*\Pi^{+*}\db{X^+}$ is a full subcategory in $\db{X}$. It is admissible, that is, there exists a non--zero left orthogonal to it.

Let $Z\in{}^{\perp}\Pi_*\Pi^{+*}\db{X^+}$ be a non--zero object. Then
$$
{\h Z, {\Pi_*\Pi^{+*}A}}=0\quad\mbox{for any }\:A\in\db{X^+}.
$$
It follows that
$$
{\h {\Pi^{+*}A} , {\Pi^*Z\otimes\omega_{\tilde X}[3]}}^* \cong
{\h \Pi^*Z, {\Pi^{+*}A}}\cong{\h Z, {\Pi_*\Pi^{+*}A}}=0,
$$
i.e. $\Pi^*Z\otimes\omega_{\tilde X}\in{\cal A}_0^+=\langle{\cal D}_0,\;{\cal B}_0\rangle$.
Let $K\in {\cal D}_0$, then $K=K'\otimes\omega_{\tilde X}$ with $K'\in{\cal C}_0$. Further, 
$$
{\h \Pi^*Z\otimes\omega_{\tilde X}, K}\cong {\h \Pi^*Z, {K'}}=0,
$$
as $K'\in{\cal C}_0\subset D({X})^{\perp}$.

It follows that $\Pi^*Z\otimes\omega_{\tilde X}\in{\cal B}_0\subset D(X)^{\perp}$, in other words, for every $M\in\db{X}$
$$
{\h \Pi^*M, {\Pi^*Z\otimes\omega_{\tilde X}}}=0.
$$
From the other hand, by Serre duality one has:
$$
{\h \Pi^*Z[-3], {\Pi^*Z\otimes\omega_{\tilde X}}}\cong{\h \Pi^*Z, {\Pi^*Z}}^*\ne 0.
$$
This proves that  ${}^{\perp}\Pi_*\Pi^{+*}\db{X^+}$ is zero.
Therefore $\Pi_*\Pi^{+*}$is an equivalence of categories.

\bigskip

{\sc Remark.} The theorem on existence  of flip is valid only in the category of Moishezon varieties. Though we have considered here only algebraic varieties, all the same works with minor changes in the Moishezon case.


{\section{Reconstruction of a variety from the derived category of coherent sheaves.}}

We have seen above that there exist examples of different varieties having equivalent the derived categories of coherent sheaves. Does it mean that $\db{X}$ is a weak invariant of a variety? In this chapter we are going to show that this is not the case.

Specifically, we prove that a variety is uniquely determined by its category, if its anticanonical (Fano case) or canonical (general type case) class is ample. In fact, the proof indicates that obstructions to the reconstruction mostly due to (partial) triviality of the canonical class. The idea is that for good, in the above sense, varieties we can recognize the one--dimensional skyskraper sheaves in $\db{X}$, using nothing but the triangulated structure of the category. The main tool for this is the Serre functor (see ch.2).

Let $D$ be a k--linear triangulated category. Denote by $F_D$ the Serre functor in $D$ (in case it exists). Recall, that if $D=\db{X}$, where $X$ is an algebraic variety of dimension $n$, then by Serre--Grothendieck duality:
\begin{equation}\label{Ser}
F_{\db{X}}(\cdot)= (\cdot)\otimes \omega_X[n],
\end{equation}
where $\omega_X$ is the canonical sheaf on $X$.
\th{Definition}\label{po}
An object $P\in D$ is called {\sf point object} of codimention $s$, if
$$
\begin{array}{ll}
i)& F_D(P)\simeq P[s],\\
ii)& {\rm Hom}^{<0}(P\:,\; P)=0,\\
iii)& {\rm Hom}^{0}(P\:,\; P)=k.
\end{array}
$$
\eth
\th{Proposition}\label{rp}
Let $X$ be a smooth algebraic variety of dimension $n$ with the ample canonical or anticanonical sheaf. Then an object $P\in \db{X}$ is a point object, iff $P\cong{\o x}[r]$ is isomorphic (up to translation) to a one-dimensional skyscraper sheaf of a closed point $x\in X$.
\eth
{\sc Remark.} Since $X$ has an ample invertible sheaf it
is projective.
\pr
Any one--dimensional skyscraper sheaf obviously satisfies properties of a point object of the same codimension as the dimension of the variety.

Suppose now that for some object $P\in\db{X}$ properties i)--iii) of definition \ref{po} are verified.

Let ${\cal H}^i$ are cohomology sheaves of $P$. It immediately follows from i) that $s=n$ and ${\cal H}^i\otimes \omega_X={\cal H}^i$. Since $\omega_X$ is an ample or antiample sheaf, we conclude that ${\cal H}^i$ are finite
length sheaves, i.e. their support are isolated points. Sheaves with the support in different points are homologically orthogonal, therefore any  such object decomposes into direct some of those having the support of all cohomology sheaves in a single point. By iii) the object $P$ is indecomposable. Now consider the spectral sequence, which calculates ${\rm Hom}^m( P\:, \; P)$ by ${\rm Ext}^i( {\cal H}^j\:,\; {\cal H}^k)$:
$$
E^{p,q}_2= \bigoplus_{k-j=q}{\rm Ext}^p({\cal H}^j\:,\; {\cal H}^k) \Longrightarrow {\rm Hom}^m( P\:, \; P).
$$

Let us mention that for any two finite length sheaves having the same single point as their support, there exists a  non--trivial homomorphism from one to the other, which 
sends   generators of the first one to the socle of the second.

Considering ${\rm Hom}^m({\cal H}^j \:, \; {\cal H}^k)$ with minimal $k-j$, we observe that this non--trivial space survives at $E_{\infty}$, hence by ii) $k-j=0$. That means that all but one cohomology sheaves are zero. Moreover, iii) implies that this sheaf is a one--dimensional skyscraper. This concludes the proof.

Now having the skyscrapers we are able to reconstruct the invertible sheaves.
\th{Definition}\label{inv}
An {\sf object} $l\in D$ is called {\sf invertible} if for any point object $P\in D$ there exists $s\in {\bf Z}$ such that
$$
\begin{array}{lll}
i)& {\rm Hom}^s(L\:,\; P)=k,&\\
ii)& {\rm Hom}^i(L\:,\; P)=0, & \quad\mbox{for}\;i\ne s\\
\end{array}
$$
\eth
\th{Proposition}\label{rin}
Let $X$ be a smooth irreducible algebraic variety. Assume that all point objects have the form ${\o x}[s]$ for some $x\in X, s\in {\bf Z}$. Then an object $L\in D$ is invertible, iff $L\cong {\cal L}[p]$ for some linear vector bundle ${\cal L}$ on $X$.
\eth
\pr
For a linear bundle ${\cal L}$ we have:
$$
{\rm Hom}({\cal L}\:, \;{\o x})=k, \quad{\rm Ext}^i({\cal L}\:, \; {\o x})=0 \;\mbox{if }\:i\ne 0.
$$
Therefore, if $L={\cal L}[s]$, then it is an invertible object.

Now let ${\cal H}^i$ are the cohomology sheaves for an invertible object ${\cal L}$. Consider the spectral sequence, which calculates ${\rm Hom}^.({\cal L}\:, \;{\o x})$ for a point $x\in X$ by means of ${\rm Hom}^i({\cal H}^j\:, \;{\o x})$:
$$
{\rm Hom}^p({\cal H}^q\:, \;{\o x})\Longrightarrow{\rm Ext}^{p-q}({\cal L}\:, \;{\o x}).
$$

Let ${\cal H}^{q_0}$ be the non--zero cohomology sheaf with maximal index. Then for any point $x\in X$ from the support of ${\cal H}^{q_0}$ ${\rm Hom}({\cal H}^{q_0}\:, \;{\o x})\ne 0$. But both  ${\rm Hom}({\cal H}^{q_0}\:, \;{\o x})$ and ${\rm Ext}^1({\cal H}^{q_0}\:, \;{\o x})$ are intact by differential of the spectral sequence. Therefore, by definition of an invertible object we obtain that for any point $x$ from the support of ${\cal H}^{q_0}$ 
$$
\begin{array}{ll}
a)& {\rm Hom}({\cal H}^{q_0}\:,\;{\o x})=k,\\
b)& {\rm Ext}^1({\cal H}^{q_0}\:,\;{\o x})=0.
\end{array}
$$

Since $X$ is smooth and connected it follows the ${\cal H}^{q_0}$ is a locally free one dimensional sheaf on $X$.

This implies that ${\rm Ext}^i({\cal H}^{q_0}\:,\;{\o x})=0$  for $i>0$ and ${\rm Hom}({\cal H}^{q_0-1}\:,\;{\o x})$ are intact by differentials of the spectral sequence. This means that ${\rm Hom}({\cal H}^{q_0-1}\:,\;{\o x})=0$, for any $x\in X$, i.e. ${\cal H}^{q_0-1}=0$. Repeting this argument for ${\cal H}^q$ with smaller $q$, we easily see that all ${\cal H}^q$, except $q=q_0$, are zero. This proves the proposition.
  
Now we are ready to prove the reconstruction theorem. Linear bundles help us to `glue' points together.
\th{Theorem}\label{rec}
Let $X$ be a smooth irreducible projective variety with ample canonical or anticanonical sheaf. If $D=\db{X}$ is equivalent
  as a triangulated category to $\db{X'}$ for some other smooth algebraic variety $X'$, then $X$ is isomorphic to $X'$.
\eth
This theorem is stronger than just a reconstruction for a variety with ample canonical or anticanonical sheaf from its derived category.

One have to be careful: since $X'$ might not have ample canonical or anticanonical sheaf, the situation is not symmetric with respect to $X$ and $X'$.

We divide the proof in several steps, so that the reconstruction procedure was transparent.
\pr
Step 1. Denote ${\cal P}_D$ the set of isomorphism classes of the point objects in $D$, ${\cal P}_X$ the set of isomorphism  classes of objects in $\db{X}$
$$
{\cal P}_X:=\Bigl\{ {\o x}[k]\;\Bigl |\: x\in X, k\in {\bf Z}\Bigl\}.
$$
By proposition \ref{rp} ${\cal P}_D\cong{\cal P}_X$. Obviously, ${\cal P}_X'\subset {\cal P}_D$. Suppose that there is an object $P\subset {\cal P}_D$, which is not contained in ${\cal P}_X'$. Since ${\cal P}_D\cong{\cal P}_X$, any two objects in ${\cal P}_D$ either are homologically mutually orthogonal or belong to a common orbit with respect to the translation functor. It follows that $P\in\db{X'}$ is orthogonal to any skyscraper sheaf
 ${\cal O}_{x'}, x'\in X'$. Hence $P$ is zero. Therefore, ${\cal P}_X'=
{\cal P}_D= {\cal P}_X$.

Step 2.  Denote by ${\cal L}_D$ the set of isomorphism classes of invertible objects in $D$, ${\cal L}_X$ the set of
isomorphism classes of objects in $\db{X}$
$$
{\cal L}_X:=\Bigl\{ L[k]\;\Bigl |\: L\:\mbox{is linear bundle on } X , k\in {\bf Z}\Bigl\}. 
$$
By step 1 both varieties $X$ and $X'$ satisfy the assumptions of proposition \ref{rin}. It follows that ${\cal L}_X= {\cal L}_D= {\cal L}_X'$.

Step 3.  Let us fix some invertible object $L_0$ in $D$, which is a linear bundle on $X$. By step 2 $L_0$ can be considered, up to translation, as a linear bundle on $X'$. Moreover, changing if necessary, the equivalence $\db{X}\simeq\db{X'}$, by the translation functor, we can assume that $L_0$, regarded as an object on $X'$, is a genuine linear bundle.

Obviously, by step 1 the set $p_D\subset P_D$
$$
p_D:=\Bigl\{ P\in P_D \;\Bigl |\: {\h L, P}=k\Bigl\}
$$
coincides with both sets $p_X=\{{\o x}, x\in X\}$ and $p_X'=\{{\cal O}_{x'}, x'\in X'\}$. This gives us a pointwise identification of $X$ with $X'$.

Step 4.  Let now $l_X$ (resp., $l_X'$) be the subset in ${\cal L}_D$ of linear bundles on $X$ (resp., on $X'$).

They can be recognized from the triangulated structure of $D$ as follows:
$$
l_X'=l_X=l_D:=\Bigl\{ L\in {\cal L}\;\Bigl |\:{\h L, P}=k \:\mbox{for any } P\in p_D \Bigl\}.
$$
For $\alpha\in{\h L_1, {L_2}}$, where $L_1, L_2\in l_D$, and $P\in p_D$, denote by $\alpha^*_P$ the induced morphism:
$$
\alpha^*_P: {\h L_2, P}\lto {\h L_1, P},
$$
and by $U_{\alpha}$ the subset of those objects $P$ from $p_D$ for which $\alpha^*_P\ne 0$. By \cite{Il} any algebraic variety has an ample system of linear bundles. This means that $U_{\alpha}$, for all ${\alpha}, L_1, L_2$, gives a base for the Zariski topologies on both $X$ and $X'$. This means that the topologies on $X$ and $X'$ coincide.

Step 5.  Since codimension of all point objects are equal to dimension of $X$ and of $X'$, we have $dimX=dimX'$. Then, formula (\ref{Ser}) for the Serre functor shows that the operations of twisting by the canonical sheaf of $X$ and $X'$ induce equal transformations on the set $l_D$.

Let $L_i=F^iL_0[-ni]$. Then $\{L_i\}$ is the orbit of $L_0$ with respect to twisting by the canonical sheaf.

Since $\omega_X$ is ample or antiample, the set of all $U_{\alpha}$, where ${\alpha}$ runs all elements from ${\h L_i, {L_j}}, i,j\in {\bf Z}$, is the base of the Zariski topology on $X$, hence, by step 4, on $X'$. That means that canonical sheaf of $X'$ is also ample or, respectively, antiample.

This means that if we consider the graded algebra $A$ with graded components
$$
A_i={\h L_0, {L_i}}
$$
and with obvious ring structure,
then ${\bf Proj}A=X=X'$. This finishes the proof.

The problem of reconstructing of a variety from its derived category is closely related to the problem of computing the
group of auto-equivalences for $\db{X}$. For ample canonical or anticanonical class we have the following
\th{Theorem}\label{aut}
Let $X$ be a smooth irreducible projective variety with
ample canonical or anticanonical class. Then the group
of isomorphism classes of exact auto-equivalences $\db{X} \to
 \db{X} $ is generated by the automorphisms of the variety, the twists by linear bundles and the translations.
\eth

\pr
Assume for definiteness that the canonical class is ample.
Let us look more carefully at the proof of theorem \ref{rec}
for the case $X=X'$. At step 3 we can choose $L_{0}=\cal{O}$,
and using twists by linear bundles and translations, we can
assume that our functor takes  $\cal{O}$ to $\cal{O}$. Then, 
step 5 gives us an automorphism of the canonical ring. Since the canonical class is ample, automorphisms of the ring are in one-to-one correspondence with those of the variety.
Therefore, composing  our functor with an automorphism of the
 variety, we can assume that it induces the trivial automorphism
of the canonical ring. 

Thus we have a functor, which takes the
trivial linear bundle and any power of the canonical bundle to themselves and preserves homomorphisms between all these
bundles. Such a functor is isomorphic to the identity functor.

 Indeed, it preserves the abelian subcategory
of pure sheaves, because the sheaves can be characterized as
the objects having trivial ${\rm Hom}^{i}$, for $i\ne 0$, from  a sufficiently
negative power of the canonical sheaf. Any sheaf has a resolution by direct sums of powers of the canonical class.
Our functor takes such a resolution to  isomorphic one,
i.e. any sheaf goes to  isomorphic one. Since the sheaves
generate the derived category, we are done.

The problem of computing the group of auto-equivalences
for the case of non-ample canonical or anticanonical class
seems to be of considerable interest.

\end{document}